\documentclass[12pt,reqno]{amsart}

\usepackage{amssymb,bbm,esint,units,url}
\usepackage{booktabs}
\usepackage[bookmarks=false,unicode,colorlinks,urlcolor=red,citecolor=red,linkcolor=blue]{hyperref}
\usepackage{aliascnt}
\usepackage{doi}
\usepackage{tikz}
\usepackage{pgfplots} 
\usepackage{xfrac}

\usepackage{graphicx,color}
\usepackage{graphics}

\definecolor{dgreen}{rgb}{0,0.5,0.1}

\usepackage[normalem]{ulem}
\usepackage[margin=1.25in]{geometry}

\newtheorem{theorem}{Theorem}[section]

\newaliascnt{corx}{thmx}

\aliascntresetthe{corx}

\newaliascnt{lemma}{theorem}
\newtheorem{lemma}[lemma]{Lemma}
\aliascntresetthe{lemma}

\newaliascnt{proposition}{theorem}
\newtheorem{proposition}[proposition]{Proposition}
\aliascntresetthe{proposition}

\newaliascnt{corollary}{theorem}

\aliascntresetthe{corollary}

\newaliascnt{conjecture}{theorem}

\aliascntresetthe{conjecture}

\newaliascnt{example}{theorem}

\aliascntresetthe{example}

\newaliascnt{question}{theorem}

\aliascntresetthe{question}

\makeatletter
\def\tagform@#1{\maketag@@@{\ignorespaces#1\unskip\@@italiccorr}}
\let\orgtheequation\theequation
\def\theequation{(\orgtheequation)}
\makeatother
\def\equationautorefname~{}
\newcommand{\arxiv}[1]{%
\href{https://arxiv.org/abs/#1}{ArXiv:#1}}

\newcommand{\BS}{{\textup{BS}}}
\newcommand{\FL}{{\textup{FL}}}

\newcommand{\ctK}{{\operatorname{ct}_K\,}}

\newcommand{\e}{\varepsilon}
\newcommand{\HH}{{\mathbb H}}

\newcommand{\LL}{{\mathcal L}}

\newcommand{\R}{{\mathbb R}}

\renewcommand{\SS}{{\mathbb S}}
\newcommand{\sign}{\operatorname{sign}}
\newcommand{\snK}{{\operatorname{sn}_K\,}}

\newcommand{\I}{I}
\newcommand{\II}{{I{\!}I}}
\newcommand{\III}{{I{\!}I{\!}I}}
\newcommand{\IV}{{I{\!}V}}

\newcommand\soutb{\bgroup\markoverwith{\textcolor{blue}{\rule[.5ex]{2pt}{1pt}}}\ULon}

\begin{document}

\title[Maximizing the second Robin eigenvalue]{Maximizing the second Robin eigenvalue of simply connected curved membranes}
\author[]{Jeffrey J. Langford and Richard S. Laugesen}
\address{Department of Mathematics, Bucknell University, Lewisburg, PA 17837, U.S.A.}
\email{jjl026@bucknell.edu}
\address{Department of Mathematics, University of Illinois, Urbana,
IL 61801, U.S.A.}
\email{Laugesen@illinois.edu}
\date{\today}

\keywords{vibrating membrane, curvature bound, isoperimetric inequality, Laplacian eigenfunction, Laplace--Beltrami, simply connected surface, spherical cap, hyperbolic disk}
\subjclass[2010]{\text{Primary 35P15. Secondary 58J50}}

\begin{abstract}
The second eigenvalue of the Robin Laplacian is shown to be maximal for a spherical cap among simply connected Jordan domains on the $2$-sphere, for substantial intervals of positive and negative Robin parameters and areas. Geodesic disks in the hyperbolic plane similarly maximize the eigenvalue on a natural interval of negative Robin parameters. These theorems extend work of Freitas and Laugesen from the Euclidean case (zero curvature) and the authors' hyperbolic and spherical results for Neumann eigenvalues (zero Robin parameter). 

Complicating the picture is the numerically observed fact that the second Robin eigenfunction on a large spherical cap is purely radial, with no angular dependence, when the Robin parameter lies in a certain negative interval depending on the cap aperture. 
\end{abstract}

\maketitle

{\noindent \it Dedicated to the memory of my friend and mentor Peter Duren, who generously shared his knowledge of and fondness for special functions and conformal mappings. -- R.S.L.}

\section{\bf Introduction\label{intro}}  

%\subsection*{Second Robin eigenvalue}
Does the spherical cap maximize the second tone of vibration among membranes of given area on the sphere, subject to elastic boundary constraints? To formulate the problem mathematically, consider the second eigenvalue of the Laplacian under Robin boundary conditions on a spherical domain of given area. We show for a substantial range of areas and Robin parameters that the second eigenvalue is largest when the domain is a spherical cap. 

The analogous Euclidean result was proved by Freitas and Laugesen \cite{FL20,FL21}, building on Neumann techniques of Szeg\H{o} \cite{S54} and Weinberger \cite{W56}. 

The spherical situation in this paper is more difficult because the second Robin eigenfunction need not have angular dependence --- it can be purely radial. When the eigenfunction does have angular dependence, its radial part need not be monotonic: it can increase and then decrease and then increase once again. We handle such complications by building on our proof for the second spherical Neumann eigenvalue \cite{LL22a}, where we showed that the spherical cap is the maximizer among simply connected domains on the $2$-sphere of given area provided the domain covers less than $16/17 \simeq 94$\% of the whole sphere. That Neumann theorem improved on the 50\% result of Bandle \cite{B72,B80}, and thus required techniques applicable to caps beyond the hemisphere.  

\subsection*{The main theorem} 
A \emph{Jordan--Lipschitz surface} is a simply-connected, bounded planar domain $\Omega$ with Lipschitz boundary that is a Jordan curve, endowed with a mass density or weight $\omega \in C^2(\Omega) \cap C(\overline{\Omega})$ that is positive on $\overline{\Omega}$. We write $\Omega_\omega$ when it is desirable to indicate the weight. The weight generates a metric $\omega \, |dz|^2$ with area
\[
|\Omega_\omega| = \int_\Omega \omega \, dA 
\] 
and boundary length
\[
L_\omega = \int_{\partial \Omega} \sqrt{\omega} \, ds .
\]
The curvature of the surface is less than or equal to a constant $K$ if  
\[
- \frac{\Delta \log \omega}{2\omega} \leq K .
\]
For the Laplace--Beltrami operator $\omega^{-1} \Delta$, the Robin eigenvalue problem is
\[
\left\{ 
\begin{aligned}
- \Delta u & = \lambda \omega u \qquad \text{in $\Omega$,} \\
- \frac{\partial u}{\partial n} & = \alpha \sqrt{\omega} u \quad \text{on $\partial \Omega$,} 
\end{aligned}
\right.
\]
where $\partial/\partial n$ is the Euclidean normal derivative in the outward direction and $\alpha \in \R$ is the Robin parameter. The eigenvalues satisfy 
\[
\lambda_1(\Omega_\omega,\alpha) < \lambda_2(\Omega_\omega,\alpha) \leq \lambda_3(\Omega_\omega,\alpha) \leq \dots \to \infty , 
\]
with variational characterization  
\begin{equation} \label{eq:rayleigh}
\lambda_k(\Omega_\omega,\alpha) = \min_\LL \max_{u \in \LL \setminus \{ 0 \}} \frac{\int_\Omega |\nabla u|^2 \, dA + \alpha \int_{\partial \Omega} u^2 \, \sqrt{\omega} \, ds}{\int_\Omega u^2 \, \omega \, dA} 
\end{equation}
where $\LL$ ranges over $k$-dimensional subspaces of $W^{1,2}(\Omega)$. The Sobolev space imbeds compactly into $L^2(\omega \, dA)$ by the Lipschitz assumption, justifying discrete spectrum. 

This paper aims to maximize the second eigenvalue $\lambda_2$. Write $M_K$ for the complete $2$-dimensional surface of constant curvature $K$, so that $M_K$  can be identified with a sphere when $K>0$, the Euclidean plane when $K=0$, and a hyperbolic or Poincar\'{e} disk when $K<0$. Their Laplace--Beltrami operators are recalled in \autoref{sec:coordinates}. 

\begin{figure}
\begin{center}
\includegraphics[scale=0.6]{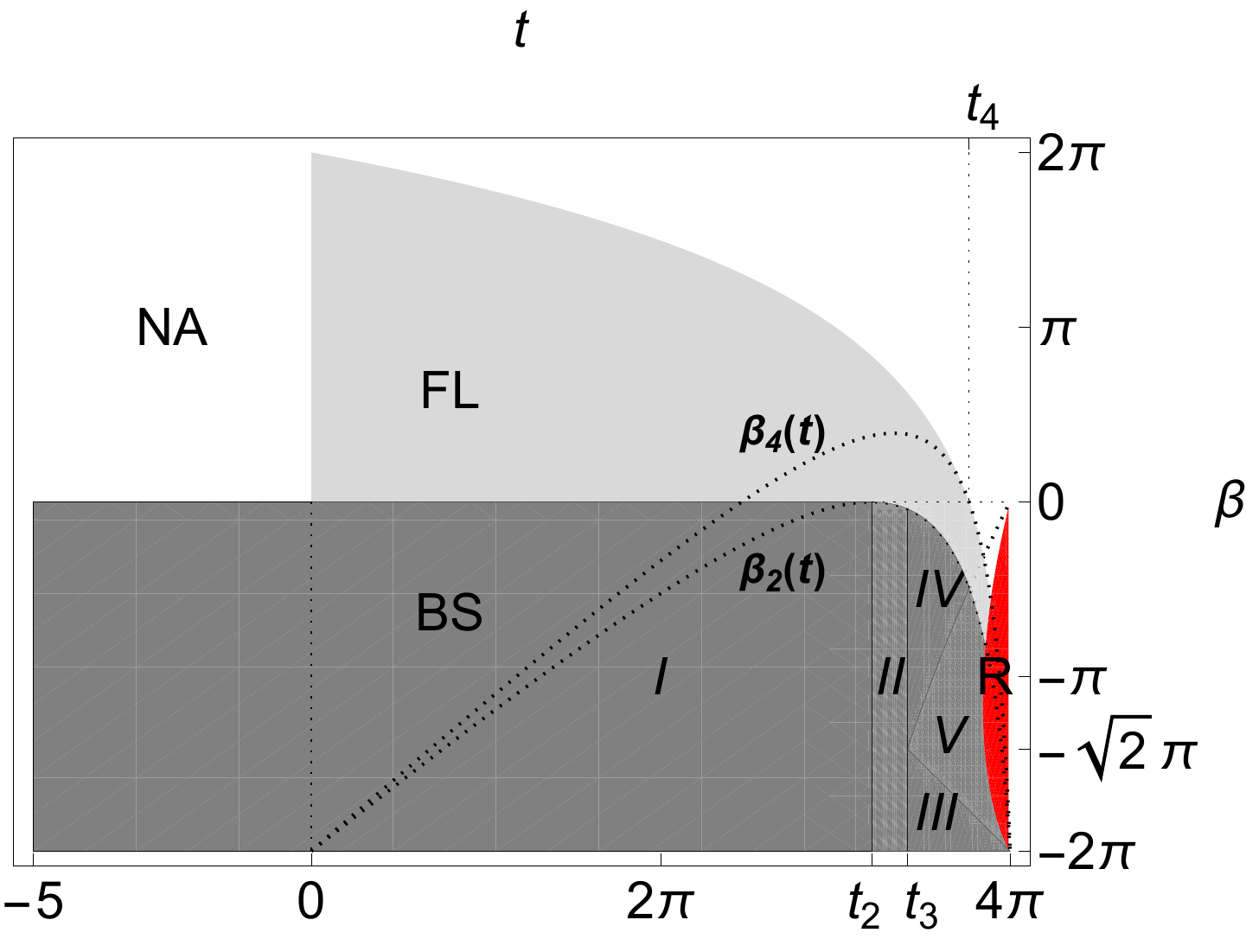}
\caption{
\label{fig:BS-FL} 
\autoref{th:robin} proves that the geodesic disk maximizes the second Robin eigenvalue for parameters lying in $\BS \cup \FL$. Sets $I$--$V$ lie in $\BS$ by \autoref{th:bs}, although the plotting of $V$ requires some numerical work. The conditions for $\FL$ are verified numerically too, as explained in \autoref{sec:construction}. \emph{Notes.} Set $I$ is the halfstrip $(-\infty,t_2] \times [-2\pi,0]$. The upper left boundary of $V$ is the graph of $t(4\pi-t)/(2\pi-t)$, which hits the rightmost point of set $\IV$ at coordinates $t \simeq 11.841$ and $\beta \simeq -1.544$, before the graph continues to increase to the point $(4\pi,0)$. The upper boundary of $\FL$ intersects the horizontal axis at $(t_4,0)$. Our method is not applicable in the second quadrant, labelled NA, or in the red region labelled R on the far right side, where the second eigenfunction is found numerically to be purely radial (\autoref{fig:angular} below).
} 
\end{center}
\end{figure}

\autoref{th:robin} below says that a constant curvature disk maximizes the second Robin eigenvalue if the curvature of the surface is bounded above and the area and the Robin parameter lie in certain regions of parameter space: a ``Bandle--Szeg\H{o}'' set $\BS$ and a ``front-loaded'' set $\FL$. These two-dimensional parameter regions are specified precisely in \autoref{sec:regions} and illustrated in \autoref{fig:BS-FL}. They involve the horizontal coordinates:
\begin{align*}
t_2 & = 4\pi \sin^2 \Theta_2/2 \simeq 10.081 \qquad \text{(defined in \autoref{sec:regions}),} \\
t_3 & = 4\pi \sin^2 3\pi/8 = (2+\sqrt{2})\pi \simeq 10.726 , \\
t_4 & \simeq 11.828 \qquad \qquad \text{(approx.\ $(16/17) 4\pi$, defined in \autoref{sec:construction}).}
\end{align*}
The $\BS$ and $\FL$ sets lie to the left of $t=4\pi$ and above $\beta=-2\pi$, so that the next theorem implicitly imposes an area restriction $|\Omega_\omega| K < 4\pi$ and parameter restriction $\beta/L_\omega \geq -2\pi$. 
\begin{theorem}[Second Robin eigenvalue is maximal for constant curvature disk] \label{th:robin}
Assume $K \in \R$ and $\Omega_\omega$ is a Jordan--Lipschitz surface with curvature $\leq K$. If $(|\Omega_\omega| K,\beta) \in \BS \cup \FL$ then
\[
\lambda_2(\Omega_\omega,\beta/L_\omega) \leq \lambda_2(D_K,\beta/L_K) 
\]
where $D_K$ is a geodesic disk in the constant curvature space $M_K$ whose boundary length is denoted $L_K$ and whose area is chosen to equal $|\Omega_\omega|$. If in addition $\beta>-2\pi$, then equality holds if and only if $\Omega_\omega$ is isometric to the constant curvature disk $D_K$.
\end{theorem}
Scaling the Robin parameter in the theorem by boundary length with $\alpha=\beta/L_\omega$ makes a natural choice, since the  parameter $\alpha$ in the Robin boundary condition must have dimension matching that of the normal derivative $\partial/\partial n$, namely $1/\text{length}$. 

\begin{figure}
\begin{center}
\includegraphics[scale=0.8]{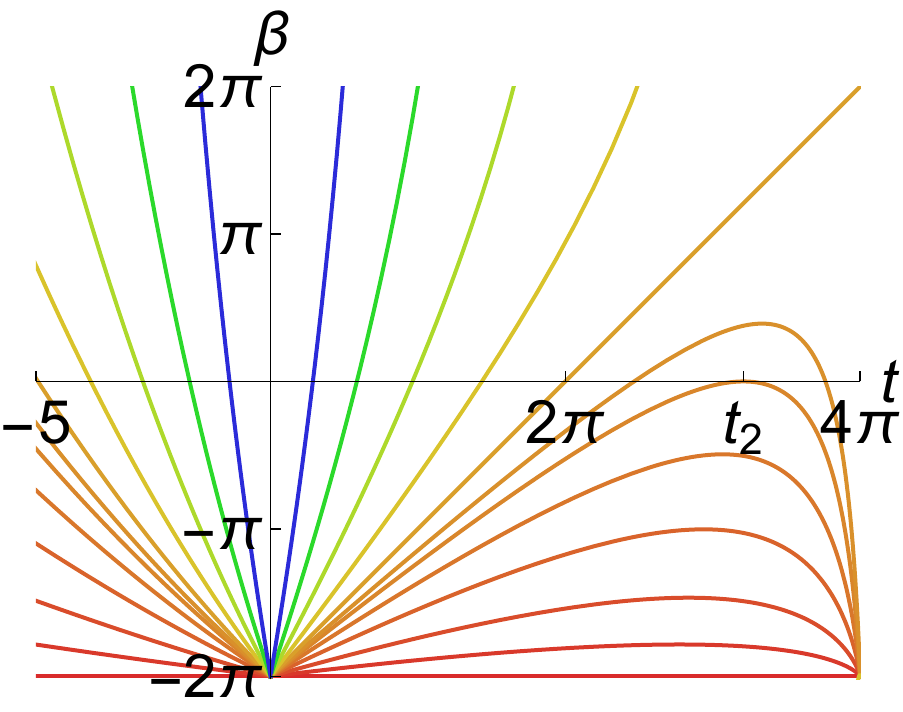}
\caption{
\label{fig:contourplot} 
Contour plot of the lowest ``angular'' eigenvalue of a geodesic disk of radius $\Theta$, whose eigenfunction has the form $g(\theta) \cos \phi$. Positive $t$: unit sphere, area of spherical cap is $t=4\pi \sin^2 \Theta/2$, Robin parameter $\alpha=\beta/2\pi \sin \Theta$. Negative $t$: hyperbolic plane, area of geodesic disk is $-t=4\pi \sinh^2 \Theta/2$, Robin parameter $\alpha=\beta/2\pi \sinh \Theta$. The eigenvalue increases as one moves upward in the figure. The horizontal line at height $-2\pi$ is the contour $\lambda=0$, corresponding to a Steklov eigenvalue for the disk. The straight line through $(2\pi,0)$ is the contour with eigenvalue $2$, corresponding to the eigenfunction $u=\sin \theta \cos \phi$ (the first spherical harmonic, that is, the coordinate function $x_1$ in $\R^3$). The contour touching the horizontal axis at $t_2$ is the graph $\beta_2(t)$ plotted in the previous figure.  For details on the construction see \autoref{sec:construction}. 
}
\end{center}
\end{figure}

The proof is in \autoref{sec:robin2proof}. On the sets $\BS$ and $\FL$, hypothesis \eqref{cond1} below ensures that the second eigenvalue $\lambda_2$ of the geodesic disk is the lowest ``angular'' eigenvalue. Thus the upper bound in the theorem is computable by separation of variables using roots of associated Legendre functions (\autoref{sec:Legendre}). Level sets of the lowest angular eigenvalue are shown in \autoref{fig:contourplot}, as a function of $\beta$ and the signed area $t=|D_K|K$. 

\subsection*{Open problems}

\subsubsection*{Problem 1 --- spherical} In the first quadrant of \autoref{fig:BS-FL}, for domains on the sphere with positive Robin parameter, can a larger region be found on which \autoref{th:robin} holds? The $\FL$ region gives a sufficient condition but is presumably not necessary. We conjecture the theorem should hold on some larger region that attains a vertex at $(4\pi,0)$. If true, then the theorem would apply in particular to the second Neumann eigenvalue on all simply connected spherical domains, with no restriction on the area. We raised that Neumann conjecture earlier \cite[Conjecture 1.2]{LL22a}. 

Does the cap maximize the second eigenvalue in the exceptional region of the fourth quadrant in \autoref{fig:BS-FL}, that is, for spherical domains with negative Robin parameter whose second eigenfunction is radial?  

\subsubsection*{Problem 2 --- hyperbolic} Our theorem does not apply in the second quadrant, for domains in hyperbolic space with positive Robin parameter. The obstacle resides in the ratio-of-areas \autoref{lem:B/A}, which holds only with nonnegative curvature. Surely the theorem itself continues to hold for a large part of the second quadrant? 

\subsection*{Prior work maximizing second eigenvalue for simply connected domains}
The original result of Szeg\H{o} \cite{S54} corresponds to the origin in \autoref{fig:BS-FL}, since he handled simply connected Euclidean domains ($K=0$) with Neumann boundary condition ($\beta=0$). Later developments by Bandle \cite{B72,B80} correspond to the interval $(-\infty,2\pi]$ on the horizontal axis in \autoref{fig:BS-FL}, that is, surfaces with Neumann boundary condition and curvature bounded above by $K$ and satisfying $|\Omega_\omega| K \leq 2\pi$. Most recently, our paper \cite{LL22a} extended Bandle's Neumann theorem to the larger interval $(-\infty,t_4]$ on the horizontal axis. The Freitas--Laugesen paper \cite{FL20} handled the interval $[-2\pi,2\pi]$ on the vertical axis in \autoref{fig:BS-FL}, in other words, it handled Euclidean domains with Robin parameter $\beta/L$ where $|\beta| \leq 2\pi$. 

The papers by Freitas--Laugesen \cite{FL20} and Langford--Laugesen \cite{LL22a} relaxed the eigenfunction monotonicity assumption that was crucial to Bandle and Szeg\H{o}'s work, by developing a modified functional that ``front loads'' the monotonicity requirement: one allows the radial part of the eigenfunction to decrease after it has first increased sufficiently. This behavior of the radial part distinguishes the two regions in \autoref{th:robin}: $\BS$ covers situations where the radial part of the eigenfunction is monotonic, and $\FL$ applies in many cases where it is not monotonic.  

In order for our methods to work, the second Robin eigenfunction of the geodesic disk must have angular dependence. That angularity requirement is built into the definitions of $\BS$ and $\FL$ in \autoref{sec:regions}. Perhaps surprisingly, the second eigenfunction can fail to have angular dependence. Numerical work shows: 
\begin{quote}
the second Robin eigenfunction on a spherical cap is purely radial (no angular dependence) when the cap fills almost the full sphere and the Robin parameter is negative and lies in a certain interval.
\end{quote}
This exceptional region of parameter space is shown in red in \autoref{fig:BS-FL}, based on the underlying plot in \autoref{fig:angular} later in the paper. 

\subsection*{Prior work maximizing the second eigenvalue for arbitrary domains} A parallel strand of research has aimed to maximize the second eigenvalue of the Laplacian for domains in all dimensions, without requiring that the domains be simply connected. Weinberger \cite{W56} showed in Euclidean space that the second Neumann eigenvalue is maximal for the ball of the same volume. The analogous result holds for subdomains of hyperbolic space by Chavel \cite{C80}, \cite[p. 94]{C84} (see also \cite{AB95,X95}), and for subdomains of the sphere that fill at most half the sphere and either contain no antipodal point-pairs (Ashbaugh and Benguria \cite[Theorem 5.1]{AB95}) or else lie outside a spherical cap of the same area (Bucur, Martinet and Nahon \cite[Corollary 3]{BMN22}). See also Wang \cite{W19} for a variable curvature result. Interestingly, the spherical cap does not always maximize the second Neumann eigenvalue among domains in $\SS^2$ that are permitted to have holes (not simply connected), as Martinet \cite{M22} has shown by numerical counterexamples for domains having large enough area. 

For the second Robin eigenvalue with a certain range of negative Robin parameters, the geodesic ball is again the maximizer among Euclidean domains by Freitas and Laugesen \cite{FL21}, whose method was extended to hyperbolic space for a smaller parameter range by Li, Wang and Wu \cite{LWW22}. 

The Bandle--Szeg\H{o} conformal mapping approach in this paper is better than the Weinberger-type approach in two key respects, for Neumann and Robin eigenvalues on simply connected subdomains of the $2$-sphere: it can treat positive Robin parameters and can handle domains with area greater than half that of the sphere. 

\subsection*{Prior work extremizing first and third Robin eigenvalues} 
To place the current paper in context, we remark that the first Robin eigenvalue too can be extremized. The sensible question now concerns minimization. The geodesic ball provides the minimizer among arbitrary domains of given volume, in spaces of constant curvature in every dimension, assuming the Robin parameter is positive. That result in Euclidean space is due to Bossel \cite{B88} and Daners \cite{D06}, and on spheres and hyperbolic space to Chen, Cheng, and Li \cite{CCL23}.

The third Robin eigenvalue is maximized by a disjoint union of disks (in a limiting sense), among simply connected planar domains, as proved by Girouard and Laugesen \cite{GL21} for a range of negative Robin parameters. The maximizer among arbitrary Euclidean domains is unknown, although numerical work does suggest it is connected \cite[Figure 4]{AFK17}. Maximizing domains in hyperbolic space or the sphere are not known. 

Maximization of the third Neumann eigenvalue (zero Robin parameter) is much better understood: the optimal shape is a disjoint union of two equal-sized geodesic balls, as shown for simply connected planar domains by Girouard, Nadirashvili and Polterovich \cite{GNP09,GP10} and by Bucur and Henrot \cite{BH19} for arbitrary Euclidean domains, and for domains in hyperbolic space by Freitas and Laugesen \cite{FL22} and on the sphere by Bucur, Martinet and Nahon \cite{BMN22}. 

An excellent survey article on Robin spectral problems can be found in the work of Bucur, Freitas and Kennedy \cite{BFK17}.

\section{\bf Laplacians on the hyperbolic space, plane and sphere} \label{sec:coordinates} 

On $2$-dimensional hyperbolic space $\HH^2$ with curvature $-1$, let $\theta$ be the geodesic distance from the origin and $\phi \in (-\pi,\pi]$ be the angle measured around the origin. In the Euclidean plane, use polar coordinates with $\theta$ being the radial variable and  $\phi$ the angle around the origin. On the unit sphere $\SS^2$ with curvature $+1$, write $\theta \in [0,\pi]$ for the angle measured from the positive $z$-axis, that is, the geodesic distance from the north pole, and write $\phi \in (-\pi,\pi]$ for the longitudinal angle. 

After defining
\begin{equation} \label{eq:snKdef}
\snK \theta  = 
\begin{cases}
\sin \theta & \text{if $K=+1$,} \\
\theta & \text{if $K=0$,} \\
\sinh \theta & \text{if $K=-1$,} 
\end{cases}
\end{equation}
the Laplace--Beltrami operators for the hyperbolic ($K=-1$), Euclidean ($K=0$) and spherical ($K=+1$) situations can be written in the unified form
\[
\Delta_K \, u = \frac{1}{\snK \theta} \frac{\partial\ }{\partial \theta} \left( \snK \theta \, \frac{\partial u}{\partial \theta} \right) + \frac{1}{(\snK \theta)^2} \, \frac{\partial^2 u}{\partial \phi^2} .
\]
We are particularly interested in eigenvalues of this operator on the geodesic disk $D_K(\Theta)$ of constant curvature $K$ and radius $\Theta>0$. That disk has area $4\pi (\snK \Theta/2)^2$. In the spherical situation ($K=+1$), the radius of the disk is restricted to $\Theta<\pi$. The Robin boundary condition with parameter $\alpha$ says $-\partial u/\partial \theta = \alpha u$ at $\theta=\Theta$.

Geodesic disks in hyperbolic space and the sphere are equivalent to Euclidean disks with weight function $w_K$ and constant curvature $K$: the stereographic change of variable can be found in \cite[Section 2]{LL22a}, and $w_K$ is stated later in \eqref{eq:weightdef}. Thus just as  the competitor surface $\Omega_\omega$ in \autoref{th:robin} is a Jordan-Lipschitz surface, so is the geodesic disk $D_K$ that provides the maximizer.

\section{\bf The $\BS$ and $\FL$ sets} \label{sec:regions} 

Here we define the $\BS$ and $\FL$ regions on which \autoref{th:robin} is valid, and develop conditions for belonging to those sets. 

Given $K=-1,0,+1$ as in the preceding section, denote by 
\[
\lambda_k(\Theta,\alpha) , \qquad k=1,2,3,\dots,
\]
the $k$-th eigenvalue of $\Delta_K$ on a geodesic disk $D_K(\Theta)$ with Robin parameter $\alpha$. 

\subsection*{Definition of the $\BS$ set}
The $\BS$ set consists of parameter values for which the second eigenfunction on  $D_K(\Theta)$ has angular dependence and monotonic radial part: 
\begin{align*}
\BS 
= \{ (4\pi \sin^2 \Theta/2 , \beta ) & : -2\pi \leq \beta \leq 0 , 0<\Theta<\pi, \text{and \eqref{cond1}--\eqref{cond2a} hold for $K=+1$} \} \\
\cup \, \{ (0 , \beta ) & : -2\pi \leq \beta \leq 0 , 0< \Theta < \infty, \text{and \eqref{cond1}--\eqref{cond2a} hold for $K=0$} \} \\
\cup \, \{ (-4\pi \sinh^2 \Theta/2 , \beta ) & : -2\pi \leq \beta \leq 0 , 0 < \Theta < \infty, \text{and \eqref{cond1}--\eqref{cond2a} hold for $K=-1$} \} . 
\end{align*}
Here the first coordinate $4\pi (\snK  \Theta/2)^2 K$ is the signed area of the disk $D_K(\Theta)$.

\smallskip
\noindent \textit{Angular condition}
\begin{equation}
\text{A second eigenfunction for eigenvalue $\lambda_2(\Theta,\beta/2\pi \, \snK \Theta)$ has the form $g(\theta) \cos \phi$.} \label{cond1} 
\end{equation}
\noindent \textit{Monotonic condition}
\begin{equation}
\text{$g$ and $g^\prime$ are positive on $(0,\Theta)$,} \label{cond2a} 
\end{equation}
except that $g^\prime$ might vanish at one point in the interval. In the Euclidean case ($K=0$), if the angular condition \eqref{cond1} holds for some $\Theta$ then by scaling invariance it holds for all $\Theta$, and similarly for the monotonic condition \eqref{cond2a}. 

\subsection*{Shape of the $\BS$ set}
To state the next theorem, which provides sufficient conditions for belonging to $\BS$, we need some special functions. Define 
\[
\Theta_2 \simeq 0.70\pi 
\]
to be the unique aperture $\Theta$ of a geodesic disk $D_{+1}(\Theta)$ on the unit sphere for which the second Neumann eigenvalue $\lambda_2(\Theta,0)$ equals $\csc^2 \Theta$; see \cite[Propositions 3.1, 4.2]{LL22a} and \cite[Theorem 1]{LL22b} for the construction of this number $\Theta_2$. The corresponding Neumann eigenfunctions on the cap of aperture $\Theta_2$ have the form $g_2(\theta) \cos \phi$ and $g_2(\theta) \sin \phi$, where as shown in \cite[Proposition 4.1(a)(c)]{LL22a}, the radial part $g_2$ has positive derivative:  $g_2^\prime(\theta)>0$ for all $\theta \in (0,\pi)$ except at $\Theta_2$, where the Neumann condition requires $g_2^\prime(\Theta_2)=0$. For no other aperture is the radial part of the Neumann eigenfunction increasing on the whole interval $(0,\pi)$.

\begin{figure}
\begin{center}
\includegraphics[scale=0.6]{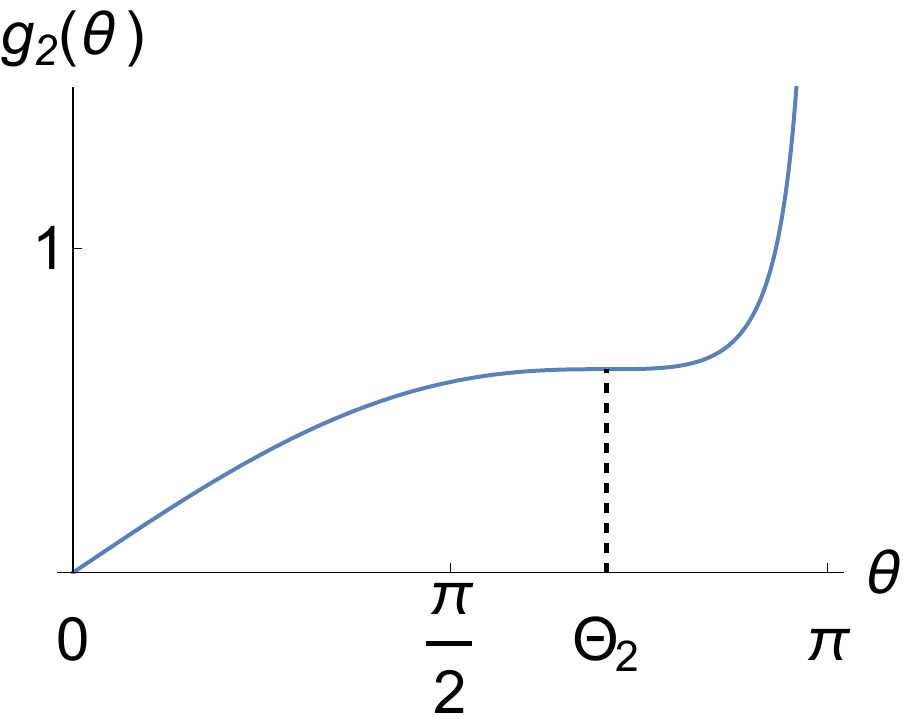}
\qquad \quad 
\includegraphics[scale=0.7]{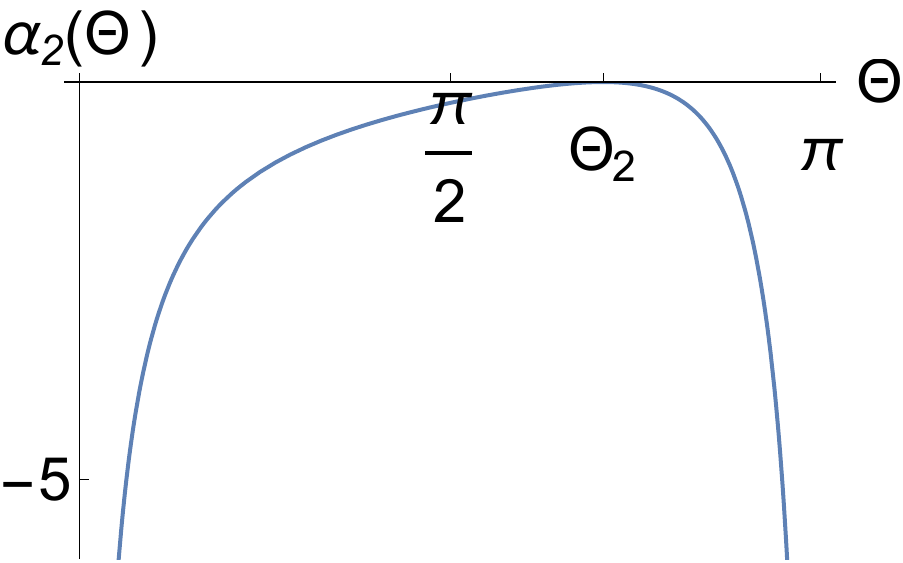}
\caption{
\label{fig:g-plot} 
\textsc{Left:} radial part $g_2(\theta)$ of the second Neumann eigenfunction for a spherical cap of aperture $\Theta_2$. \textsc{Right:} Robin parameter $\alpha_2(\Theta)=-g_2^\prime(\Theta)/g_2(\Theta)$ for $g_2$, at arbitrary aperture $\Theta$.
}
\end{center}
\end{figure}

The graph of $g_2$ on the left of \autoref{fig:g-plot} is obtained from the explicit formula (see \autoref{sec:Legendre} or \cite[Proof of Proposition 4.2]{LL22a}) that 
\[
g_2(\theta)=P_{n_2}^{-1}(\cos \theta) , \qquad \theta \in [0,\pi) ,
\]
where $P_n^m$ is the associated Legendre function and we choose $m=-1$ and take $n=n_2$ such that $g_2^\prime$ has the required property of being positive except at one point $\Theta_2$. Numerically, one finds   
\[
n_2 \simeq 0.851187 .
\]

Define 
\[
\alpha_2(\Theta) = -g_2^\prime(\Theta)/g_2(\Theta) ,
\]
so that $\alpha_2(\Theta)$ is the Robin parameter for $g_2(\theta) \cos \phi$ at the boundary of the cap of aperture $\Theta$. In particular, the Neumann condition at aperture $\Theta_2$ says $\alpha_2(\Theta_2)=0$, as seen on the right of  \autoref{fig:g-plot}. 

The area of a spherical cap of aperture $\Theta \in (0,\pi)$ is given by the strictly increasing function $t(\Theta)=4 \pi \sin^2 \Theta/2 \in (0,4\pi)$. Let   
\begin{align*}
t_2 = t(\Theta_2) & = 4\pi \sin^2 \Theta_2/2 \simeq 10.081 , \\
t_3 = t(3\pi/4) & = 4\pi \sin^2 3\pi/8 = (2+\sqrt{2})\pi \simeq 10.726 ,
\end{align*}
and define $\beta_2 : (0,4\pi) \to \R$ by  
\[
\beta_2(t(\Theta)) = (2\pi \sin \Theta) \, \alpha_2(\Theta) .
\]
Note this definition has the same form ``$\alpha=\beta/L$'' as appears in \autoref{th:robin}, since the cap of aperture $\Theta$ has boundary length $2\pi \sin \Theta$.
\begin{theorem}[Shape of the $\BS$ region] \label{th:bs}
The $\BS$ region contains the following sets:
\begin{align*}
I & = \{ (t,\beta) : -\infty < t \leq t_2 , -2\pi \leq \beta \leq 0 \} \\
\II & = \{ (t,\beta) : t_2 < t \leq t_3 , -2\pi \leq \beta \leq \beta_2(t) \} , \\
\III & = \{ (t,\beta) : t_3 < t < 4\pi , -2\pi \leq \beta \leq 2\pi - t \} , \\
\IV & = \left\{ (t,\beta) : t_3 < t < 4\pi , \, t \, \frac{4\pi-t}{2\pi-t} \leq \beta \leq \beta_2(t) \right\} , \\
V & = \left\{ (t,\beta) : t_3 < t < 4\pi , 2\pi - t \leq \beta \leq \min \left( t \, \frac{4\pi-t}{2\pi-t} , \beta_2(t) \right) \text{ and \eqref{cond1} holds} \right\} .
\end{align*}
\end{theorem}
These sets are shown in \autoref{fig:BS-FL}. The theorem is proved in \autoref{sec:bsproof}. 

\subsection*{Definition of the $\FL$ set}
The $\FL$ set comprises those parameter values for which the second eigenfunction has angular dependence and its radial part increases and then decreases in a ``front-loaded'' way with more increase than decrease, according to a certain integral criterion:
\begin{align*}
\FL 
= \{ (4\pi \sin^2 \Theta/2 , \beta ) & : -2\pi \leq \beta < \infty , 0<\Theta<\pi, \text{and \eqref{cond1}, \eqref{cond2b}, \eqref{cond3} hold for $K=+1$} \} \\
\cup \, \{ (0 , \beta ) & : -2\pi \leq \beta < \infty , 0< \Theta < \infty , \text{and \eqref{cond1}, \eqref{cond2b}, \eqref{cond3} hold for $K=0$} \} 
\end{align*}
where the angular condition \eqref{cond1} was stated above and the new conditions are as follows. 

\smallskip
\noindent \textit{Up-Down-(Up) condition}
\begin{align}
& \text{$g>0$ on $(0,\Theta]$, $g^\prime>0$ on $(0,\theta_{max})$, $g^\prime<0$ on $(\theta_{max},\theta_{min})$, $g^\prime>0$ on $(\theta_{min},\Theta)$} \label{cond2b} \\
& \hspace*{2.8cm}  \text{for some numbers $0 < \theta_{max} < \theta_{min} \leq \Theta$,} \notag 
\end{align}
(If $\theta_{min}<\Theta$ then $g$ goes up-down-up, while if $\theta_{min}=\Theta$ then the third interval $(\theta_{min},\Theta)$ is empty and $g$ goes only up-down.)

\smallskip
\noindent \textit{Front-Loaded condition}
\begin{equation}
\int_0^{\theta_{min}} g(\theta) g^\prime(\theta) \left( \snK \frac{\theta}{2} \right)^{\! \! 2} d\theta \geq 0 . \label{cond3}
\end{equation}

The $\FL$ set lies in the right halfplane, relating to the spherical case in the first and fourth quadrants and the Euclidean case on the vertical axis. The third quadrant, meaning hyperbolic with negative Robin parameter, is handled already by the $\BS$ region, thanks to the set $\I$ in \autoref{th:bs}.  

In the second quadrant, that is, for the hyperbolic case with positive Robin parameter ($t<0,\beta>0$), we can offer no result. The obstacle is that the radial part $g$ of the second eigenfunction is non-monotonic due to the positive Robin parameter, while our tool for handling non-monotonicity (\autoref{lem:B/A}) applies only to the Euclidean and spherical cases. 

\subsection*{Shape of the $\FL$ set}
The $\FL$ set extends downward from each point it contains, as seen graphically in \autoref{fig:BS-FL}. 
\begin{proposition}[Dropping downward in the $\FL$ set] \label{pr:FLdown} Suppose $-2\pi \leq \beta < \beta_*$. 

\noindent (i) Let $K=0$ and $\Theta>0$. If $(0,\beta_*) \in \FL$ then $(0,\beta) \in \BS \cup \FL$ provided the angular condition \eqref{cond1} holds for this $\beta$. 

\noindent (ii) Let $K=+1$ and $0<\Theta<\pi$. If $(4\pi \sin^2 \Theta/2,\beta_*) \in \FL$ then $(4\pi \sin^2 \Theta/2,\beta) \in \BS \cup \FL$ provided \eqref{cond1} holds for this $\beta$ and $\Theta$. 
\end{proposition}
The proposition is proved in \autoref{sec:FLdownproof}. The angular condition \eqref{cond1} holds in particular when $\beta \geq 0$, by applying \autoref{basic} later in the paper with $\alpha \geq 0$.

Some points belonging to the $\FL$ set can be established rigorously. For example, the $\BS$ set contains the line segment with $0 < t  \leq t_2$ and $\beta=0$, and the $\FL$ set contains its continuation with $t_2 < t  \leq t_4$ and $\beta=0$, by our work in the Neumann case \cite[Theorem 1.1]{LL22a}. Further, the $\FL$ set contains the vertical line segment with $t=0$ and $0 < \beta \leq 2\pi$ by a result of Freitas and Laugesen \cite[Theorem B]{FL20} in the Euclidean case. Additional first-quadrant regions in the $\FL$ set can be determined rigorously with the help of \autoref{pr:FLdown}, if desired, as explained in \autoref{sec:construction}.

\section{\bf Curvature assumptions imply area comparisons \label{curv}}

The proof of \autoref{th:robin} relies on area growth inequalities that follow from the upper curvature bound. The first inequality is due to Bandle and addresses a difference of areas. The second inequality appeared in a recent paper of ours and deals with the ratio of areas. 

The planar weight representing the sphere, Euclidean plane or hyperbolic plane is 
\begin{equation} \label{eq:weightdef}
w_K(r)=
\begin{cases}
\frac{4}{(1+r^2)^2}, & 0\leq r<\infty, \quad \text{when $K=+1$,} \\
1, & 0\leq r<\infty,\hspace{6pt} \textup{when $K=0$,} \\
\frac{4}{(1-r^2)^2}, & 0\leq r<1, \quad \text{when $K=-1$.} 
\end{cases}
\end{equation}
One checks that the curvature $- (\Delta \log w_K)/2w_K$ equals $K$ in each case. The weighted area of the Euclidean disk $D(r)$ is
\begin{align}
A(r)=|D(r)_{w_K}| &=2\pi \int_0^rw_K(s)s \,ds \notag \\ 
&=
\begin{cases}
\frac{4\pi r^2}{1+r^2}, & 0\leq r<\infty \quad \textup{when }K=+1,\\
\pi r^2, & 0\leq r<\infty,\hspace{6pt} \textup{when }K=0,\\
\frac{4\pi r^2}{1-r^2}, & 0\leq r<1, \quad \textup{when }K=-1 .
\end{cases}
\label{def:A}
\end{align}
Notice the area $A(r)$ can take any value between $0$ and $\infty$ when $K=-1,0$, and any value between $0$ and $4\pi$ when $K=+1$. 

Suppose $K=-1,0,+1$. Given a surface $\Omega$ with weight $\omega$ as in \autoref{th:robin}, choose a radius $R>0$ such that
\[
A(R)=|\Omega_\omega| ,
\]
noting in the case $K=+1$ that such an $R$ exists because the assumptions in the theorem ensure $|\Omega_\omega|<4\pi$. Take $F:D(R) \to \Omega$ to be a conformal mapping onto the simply connected domain $\Omega$. The $\omega$-weighted area of the image of the subdisk $D(r)$ is 
\begin{equation}\label{def:B}
B(r)=|F(D(r))_\omega|=\int_{F(D(r))}\omega\,dA=\int_{D(r)}(\omega \circ F) |F^\prime|^2 \, dA.
\end{equation}
Since $\Omega=F(D(R))$, we have the endpoint condition
\begin{equation}\label{eq:ABR}
B(R)=|\Omega_\omega|=A(R) .
\end{equation}
\begin{lemma}[Difference of areas; Bandle \protect{\cite[pages 44, 119]{B80}}, or see \protect{\cite[Lemma 6.1]{LL22a}}] \label{lem:A-B}
The constant curvature disk has larger area: $A(r) \geq B(r)$ for $0<r<R$. Equality holds for all $r$ if and only if $(\omega \circ F) |F^\prime|^2 \equiv w_K$.
\end{lemma}
\begin{lemma}[Ratio of areas; Langford and Laugesen \protect{\cite[Lemma 6.3]{LL22a}}] \label{lem:B/A}
If $K = 0,1$, then the area ratio $B(r)/A(r)$ is increasing. This area ratio is constant if and only if $(\omega \circ F) |F^\prime|^2 \equiv w_K$. 
\end{lemma}

\section{\bf Proof of \autoref{th:robin} --- second Robin eigenvalue maximal for constant curvature disk}
\label{sec:robin2proof}

We follow the construction of trial functions from the Neumann case by Szeg\H{o} \cite{S54} and Bandle \cite{B72,B80}. In the hyperbolic  Robin situation we can employ their method of estimating the Rayleigh quotient, under the monotonicity condition \eqref{cond2a}. The spherical situation is handled under either \eqref{cond2a} on the $\BS$ set or else the new and distinctly weaker Front-Loaded condition \eqref{cond3} on the $\FL$ set, which enables a certain integration by parts step to be adapted from \cite{FL20,LL22a}.

Without loss of generality, we may assume the upper bound $K$ on the curvature equals $-1,0$ or $+1$, since multiplying the metric by a positive constant $c$ causes the area and boundary length to change by factors of $c$ and $\sqrt{c}$, while the curvature and eigenvalue in the theorem change by $1/c$, as is clear from the Rayleigh quotient \eqref{eq:rayleigh}. 

\subsection*{Constructing trial functions} Assume $(|\Omega_\omega| K,\beta) \in \BS \cup \FL$. The constant curvature geodesic disk $D_K$ whose area equals $|\Omega_\omega|$ lies in either the hyperbolic space of curvature $K=-1$, in Euclidean space ($K=0$), or in the unit sphere (curvature $K=+1$, noting such a spherical cap $D_K$ exists since $|\Omega_\omega|<4\pi$ by hypothesis). Write $\Theta$ for the radius of that geodesic disk. Second eigenfunctions of $-\Delta_K v = \lambda v$ on $D_K=D_K(\Theta)$ with Robin parameter $\beta/2\pi \, \snK \Theta$ can be taken in the form $v_2=g(\theta) \cos \phi$ and $v_3=g(\theta) \sin \phi$ by the angular hypothesis \eqref{cond1} in the $\BS$ and $\FL$ sets, noting that since cosine gives an eigenfunction, so must sine. 

Transform the radial variable by $r=\tanh \theta/2$ (hyperbolic) or $r=\theta$ (Euclidean) or $r=\tan \theta/2$ (spherical), and similarly define $R$ in terms of $\Theta$ in each case. (In the hyperbolic case, note that $R=\tanh \Theta/2<1$.) Writing 
\[
h(r)=g(\theta) ,
\]
one calculates (see for example \cite[Section 2]{LL22a}) that the transformed eigenfunctions 
\[
f_2=h(r) \cos \phi , \qquad f_3=h(r) \sin \phi ,
\]
are second eigenfunctions of $-\Delta f = \lambda w_K f$ on the Euclidean disk $D(R)$ having the Robin parameter $\beta/2\pi R \sqrt{w_K(R)}$, where the weight $w_K$ was defined in the previous section. Here $\Delta$ is the Euclidean Laplacian. The radial part $h$ is smooth, and has $h(0)=0$ since eigenfunctions are continuous at the origin. 

This change of variable also implies that the weighted disk has the same area as the geodesic disk, $|D(R)_{w_K}|=4\pi(\snK \Theta/2)^2=|D_K(\Theta)|$, and hence by our construction has the same area as $\Omega$ with weight $\omega$, that is, $A(R)=|\Omega_\omega|$ in the notation of \autoref{curv}. Hence $A(R)=B(R)$ by formula \eqref{eq:ABR}. 

Take a conformal mapping $F:D(R) \to \Omega$ onto the simply connected domain. Conformally transplant $f_2$ and $f_3$ to $\Omega$ by letting
\[
\varphi_2=f_2\circ F^{-1} \quad \textup{and} \quad \varphi_3=f_3\circ F^{-1}.
\]
Observe $\varphi_2,\varphi_3\in H^1(\Omega)$ by boundedness of $h$ and by conformal invariance, which yields equality and finiteness of the Dirichlet integrals:  
\[
\int_{\Omega} |\nabla \varphi_2|^2 \,dA = \int_{D(R)} |\nabla f_2|^2 \,dA = \int_{D(R)} |\nabla f_3|^3 \,dA = \int_{\Omega} |\nabla \varphi_3|^2 \,dA .
\]

By a ``center of mass'' argument that goes back to Szeg\H{o} \cite{S54} (see, for example, \cite[Lemma 5]{FL20}), we may assume after precomposing the conformal map $F$ with a suitable M\"{o}bius self-map of the disk that $\varphi_2$ and $\varphi_3$ are each orthogonal in $L^2(\omega \, dA)$ to the eigenfunction $u_1$ for the eigenvalue $\lambda_1(\Omega_\omega,\beta/L_\omega)$, meaning 
\[
\int_{\Omega} \varphi_2 u_1 \, \omega\,dA = \int_{\Omega} \varphi_3 u_1 \, \omega\,dA=0 .
\]

\subsection*{Substituting into the Rayleigh quotient} Applying the variational characterization for the second eigenvalue, restricted to the space $\{ \varphi \in H^1(\Omega) : \int_{\Omega} \varphi u_1 \, \omega\,dA=0 \}$ of functions orthogonal to the first eigenfunction, one obtains using the trial functions $\varphi_2$ and $\varphi_3$ that 
\begin{equation}
\lambda_2(\Omega_\omega, \beta/L_\omega)
\leq \frac{\int_{\Omega} |\nabla \varphi_i|^2 \,dA + (\beta/L_\omega) \int_{\partial \Omega} \varphi_i^2 \sqrt{\omega} \, ds}{\int_{\Omega}\varphi_i^2 \, \omega\,dA} \qquad i=2,3. \label{eq:trialestimate}
\end{equation}
Recall here that $L_\omega=\int_\Omega \sqrt{\omega} \, ds$ is the weighted length of the boundary. 

Clearing the denominators and summing over $i=2,3$ yields that 
\[
\lambda_2(\Omega_\omega, \beta/L_\omega) 
\leq \frac{\int_{\Omega} \left( |\nabla \varphi_2|^2 + |\nabla \varphi_3|^2 \right) dA + \beta h(R)^2}{\int_{\Omega} (\varphi_2^2 +  \varphi_3^2) \, \omega\,dA} 
\]
where we used that on $\partial \Omega$, one has $\varphi_2^2+\varphi_3^2=h(R)^2$ by the definitions. Hence  
\[
\lambda_2(\Omega_\omega, \beta/L_\omega) 
\leq \frac{\int_{D(R)} \left( |\nabla f_2|^2 + |\nabla f_3|^2 \right) dA + \beta h(R)^2}{\int_{D(R)} (f_2^2 + f_3^2) \, \omega(F) |F^\prime|^2 \,dA} 
\]
by pulling the integrals back to $D(R)$ via the conformal map $F$. After substituting $f_2=h(r)\cos \phi$ and $f_3=h(r)\sin \phi$, we find
\begin{equation}
\lambda_2(\Omega_\omega, \beta/L_\omega)  
\leq \frac{\int_{D(R)} \left( h^\prime(r)^2+r^{-2}h(r)^2 \right) dA + \beta h(R)^2}{\int_{D(R)} h(r)^2 \omega(F(re^{i\phi}))|F^\prime(re^{i\phi})|^2 \, dA} . \label{eq:eigestimate1}
\end{equation}

Equality holds in \eqref{eq:trialestimate} when $\Omega=D(R), \omega=w_K, F(z)=z, \varphi_i=f_i$, and so  
\begin{equation} \label{eq:eigestimate2}
\lambda_2(D(R)_{w_K}, \beta/L_{w_K})
= \frac{\int_{D(R)} \left( h^\prime(r)^2+r^{-2}h(r)^2 \right) dA + \beta h(R)^2}{\int_{D(R)} h(r)^2 w_K(r) \, dA} .
\end{equation}

\subsection*{Nonnegativity of the numerators} The numerators in \eqref{eq:eigestimate1} and \eqref{eq:eigestimate2} are identical. We show they are positive, except in a borderline case where they equal zero. The underlying reason is that the first nonzero Steklov eigenvalue of the Euclidean disk of radius $R$ equals $1/R$; rather than relying on that interpretation, for simplicity's sake we estimate the numerator explicitly: 
\begin{align*}
\int_{D(R)} \left( h^\prime(r)^2+r^{-2}h(r)^2 \right) dA 
& \geq \int_{D(R)} 2 h^\prime(r) r^{-1}h(r) \, dA \qquad \text{since $a^2+b^2 \geq 2ab$} \\
& = 2\pi h(R)^2 \geq -\beta h(R)^2 
\end{align*}
because $h(0)=0$ and $\beta \geq -2\pi$ in the $\BS$ and $\FL$ parameter sets. Thus the numerator in \eqref{eq:eigestimate1} and \eqref{eq:eigestimate2} is nonnegative. Equivalently, the second eigenvalue of the disk is nonnegative: $\lambda_2(D(R)_{w_K}, \beta/L_{w_K}) \geq 0$.

The numerator equals zero if and only if $\beta=-2\pi$, as we now explain. Note that $h(R)=g(\Theta)>0$ by hypothesis \eqref{cond2a} for the $\BS$ set or \eqref{cond2b} for the $\FL$ set. Hence from the inequalities in the argument above we deduce that if the numerator equals zero then $\beta=-2\pi$ and $h^\prime(r) = r^{-1}h(r)$, so that $h(r)=r$ for all $r$. In the reverse direction, if $\beta=-2\pi$ then the disk has second eigenvalue $\lambda_2(D(R)_{w_K}, -2\pi/L_{w_K}) \leq 0$, since $u=r \cos \phi$ is a sign-changing eigenfunction with eigenvalue zero: $\Delta u = 0$ and at $r=R$ the Robin condition $\partial u/\partial r-(2\pi/L_{w_K})\sqrt{w_K}u=0$ holds. Hence the numerator of \eqref{eq:eigestimate2} equals $0$. 

\subsection*{Numerators positive} Suppose $\beta>-2\pi$, so that the numerators of \eqref{eq:eigestimate1} and \eqref{eq:eigestimate2} are positive. To complete the proof, it is enough to compare denominators and show
\begin{equation}\label{ineq:mainineq}
\int_{D(R)} h(r)^2 \omega(F(re^{i\phi})) |F^\prime(re^{i\phi})|^2 \, dA - \int_{D(R)} h(r)^2 w_K(r) \, dA \geq 0
\end{equation}
with equality if and only if $\omega(F(z))|F^\prime(z)|^2 \equiv w_K(z)$. (Regarding the equality statement in the theorem, notice $(\omega \circ F) |F^\prime|^2 \equiv w_K$ means the surface $\Omega$ with metric $\omega(z) |dz|^2$ is isometric via the conformal map $F$ to the disk $D(R)$ with metric $w_K(z) |dz|^2$, while in the other direction, if the two surfaces are isometric then their eigenvalues are the same and so equality holds in the theorem.) 

Recalling the definitions of the area functionals $A$ and $B$ in \eqref{def:A} and \eqref{def:B}, the left side of \eqref{ineq:mainineq} equals 
\[
-\int_0^R h(r)^2 \left( A^\prime(r)-B^\prime(r) \right) \,dr = \int_0^R 2h(r) h^\prime(r) \left( A(r)-B(r) \right)\,dr 
\]
after an integration by parts, where the boundary terms vanish because $A(0)=B(0)=0$ and also $A(R)=B(R)$ by the area normalization  \eqref{eq:ABR}. Therefore the task for \eqref{ineq:mainineq} is to show
\[
\int_0^R 2h(r) h^\prime(r) \left( A(r)-B(r) \right) dr \geq 0 
\]
with equality if and only if $(\omega \circ F)|F^\prime|^2 = w_K$. 
We convert back to the geodesic radial variable $\theta$ by substituting $h(r)=g(\theta)$, so that the goal becomes to show 
\begin{equation} \label{eq:wantnonneg}
\int_0^\Theta 2g(\theta) g^\prime(\theta) ( A - B) \, d\theta \geq 0 
\end{equation}
with equality if and only if $(\omega \circ F)|F^\prime|^2 = w_K$. Here the areas $A$ and $B$ are regarded as functions of $\theta$. Explicitly, one computes $A=4\pi (\snK \theta/2)^2$ by substituting $r=\tanh \theta/2$ or $r=\theta$ or $r=\tan \theta/2$, respectively, into the area formula \eqref{def:A}. 

\subsection*{$\BS$ case} In the hyperbolic case ($K=-1$) the hypotheses on the $\BS$ region ensure that \eqref{cond1} and \eqref{cond2a} hold, so that $g$ and $g^\prime$ are positive (except that $g^\prime$ might vanish at one point). We know $A - B \geq 0$ by \autoref{lem:A-B}, with equality for all $\theta$ if and only if $(\omega \circ F)|F^\prime|^2 = w_K$. Inequality \ref{eq:wantnonneg} and its equality statement follow immediately. The same holds in the Euclidean and spherical cases. 

\subsection*{$\FL$ case} Consider now the Euclidean and spherical cases assuming the $\FL$ set hypotheses \eqref{cond1}, \eqref{cond2b}, \eqref{cond3}. Define 
\[
G(\psi)=\int_0^\psi 2g(\theta)g'(\theta) A \,d \theta
\] 
and note $G(0)=0$. The left side of inequality \eqref{eq:wantnonneg} can be rewritten in terms of $G$ by first pulling out a factor of $A$, obtaining   
\begin{align*}
\int_0^\Theta 2g(\theta)g^\prime(\theta)A \left(1-\frac{B}{A}\right) d\theta
&=\int_0^\Theta G^\prime(\theta)\left(1-\frac{B}{A}\right) d\theta \\
&=\int_0^\Theta G(\theta) \frac{d\ }{d\theta} \left(\frac{B}{A}\right) d\theta ,
\end{align*}
where the final step uses integration by parts and the normalization that $B/A=1$ when $\theta=\Theta$, by \eqref{eq:ABR}. The area ratio is increasing, with $(d/d\theta)(B/A) \geq 0$ by \autoref{lem:B/A}, remembering that the lemma needs $K=0,+1$; equality holds for all $\theta$ if and only if $(\omega \circ F) |F^\prime|^2 \equiv w_K$. Thus  to prove \eqref{eq:wantnonneg}, it suffices to show $G(\theta) > 0$ for all $\theta \in (0,\Theta)$ except perhaps at one $\theta$ value. 

Since $G^\prime(\theta)=2g(\theta) g^\prime(\theta) A$, the Up-Down-Up hypothesis \eqref{cond2b} implies $G$ is strictly increasing for $0<\theta<\theta_{max}$, strictly decreasing for $\theta_{max}<\theta < \theta_{min}$, and strictly increasing for $\theta_{min}<\theta<\Theta$. Because $G(0)=0$ by construction and $G(\theta_{min}) \geq 0$ by the Front-Loaded hypothesis \eqref{cond3}, it follows that $G > 0$ on the interval $(0,\Theta)$ except perhaps at $\theta_{min}$, which completes the proof. 

\subsection*{Numerators equal zero} Lastly, if the numerators of \eqref{eq:eigestimate1} and \eqref{eq:eigestimate2} equal zero then those formulas imply 
\begin{equation} \label{eq:zerocase}
\lambda_2(\Omega_\omega, \beta/L_\omega) \leq 0 = \lambda_2(D(R)_{w_K}, \beta/L_{w_K}) ,
\end{equation}
so that the second eigenvalue is maximal for the constant curvature disk. Further, $h(r)=r$ for all $r$, as shown above, and $\beta=-2\pi$. 

 \subsubsection*{Remark.} The theorem asserts no equality statement when $\beta=-2\pi$. Equality obviously holds in \eqref{eq:zerocase} ``if'' $\Omega_\omega$ is isometric to the constant curvature disk $D(R)_{w_K}$, but we do not assert ``only if''. An equality statement can be developed, though, in terms of the conformal map $F$. For suppose equality holds in the eigenvalue inequality \eqref{eq:zerocase}. By imposing equality in the variational principle \eqref{eq:trialestimate}, we find the trial functions $\varphi_2$ and $\varphi_3$ must be Robin eigenfunctions on $\Omega_\omega$ with eigenvalue $0$. The weak formulation of the eigenfunction equation for $\varphi_2$, with $\beta=-2\pi$ and eigenvalue $0$, says
\[
\int_\Omega \nabla \varphi_2 \cdot \nabla \widetilde{\varphi} \, dx - \frac{2\pi}{L_\omega} \int_{\partial \Omega} \varphi_2 \widetilde{\varphi} \sqrt{\omega} \, ds = 0, \qquad \widetilde{\varphi} \in H^1(\Omega) .
\]
Adapting an argument of Freitas and Laugesen \cite[p.{\,}1038--1039]{FL20} for Jordan--Lipschitz domains, one deduces that 
\begin{equation} \label{eq:Steklovboundary}
\sqrt{(\omega \circ F)(Re^{i\phi})} |F^\prime(Re^{i\phi})| = \frac{L_\omega}{2\pi R} 
\end{equation}
for almost every $\phi \in [0,2\pi]$, and that $\log |F^\prime|$ equals the Poisson integral of its boundary values. By \eqref{eq:Steklovboundary}, the boundary function is $\log \big( L_\omega/2\pi R \sqrt{(\omega \circ F)(Re^{i\phi})} \big)$, which depends continuously on $\phi$. Hence the harmonic function $\log |F^\prime|$ extends continuously to the closed disk, and so we have proved that that when $\beta=-2\pi$ and equality holds in the theorem, the product $(\omega \circ F) |F^\prime|^2$ extends continuously to the closure of the disk and equals a constant on the boundary. But that information does not determine its values in the interior of the disk, as happened in the equality case when $\beta>-2\pi$.

\section{\bf Robin eigenfunctions on disks with constant curvature}
\label{eigenstuff}

Angular dependence of the second Robin eigenfunction on constant curvature disks is established in this section for most (but not all) parameter values, along with  monotonicity of the radial part in some (but not all) parameter regimes. 

The hyperbolic and Euclidean situations are relatively straightforward. The spherical case is not. Proofs are given later in this section and the results are subsequently applied in \autoref{sec:bsproof} to establish subsets of $\BS$, for \autoref{th:bs}. 

\subsection*{Angular dependence of the second eigenfunction: mostly but not always} Recall that $\lambda_k(\Theta,\alpha)$ is the $k$-th eigenvalue of $\Delta_K$ on a geodesic disk $D_K(\Theta)$, under a Robin boundary condition with parameter $\alpha$. In geometric terms, when $K=+1$ the geodesic disk is the spherical cap $\{ (x_1,x_2,x_3) \in \SS^2 : x_3 > \cos \Theta \}$ having aperture or geodesic radius $0<\Theta<\pi$.

\begin{proposition}[Second eigenfunction has angular dependence in most cases] \label{basic} \ \\
Let $\alpha \in \R$. If: 

(a) [Hyperbolic] $K=-1$ and $\Theta>0$, or 

(b) [Euclidean] $K=0$ and $\Theta>0$, or 

(c) [Spherical] $K=+1$, and either $0<\Theta \leq 3\pi/4$ or else  
\[
\text{$3\pi/4<\Theta<\pi$ and $\alpha \notin (\cot \Theta,\tan \Theta)$,}
\]
then the second eigenspace of $-\Delta_K \, u = \lambda u$ on the geodesic disk $D_K(\Theta)$ with Robin parameter $\alpha$ is spanned by two functions of the form 
\[
g(\theta) \cos \phi \qquad \text{and} \qquad g(\theta) \sin \phi ,
\]
where the radial part $g$ is smooth and satisfies the Robin boundary condition 
\[
g^\prime(\Theta)+\alpha g(\Theta)=0 .
\] 
\end{proposition}
The hyperbolic and Euclidean cases of the proposition include all real Robin parameters $\alpha$, and the spherical case includes all nonnegative parameters $\alpha \geq 0$, since if $3\pi/4<\Theta<\pi$ then the excluded interval from $\cot \Theta$ to $\tan \Theta$ contains only negative numbers.  

\begin{figure}
\begin{center}
\includegraphics[scale=0.6]{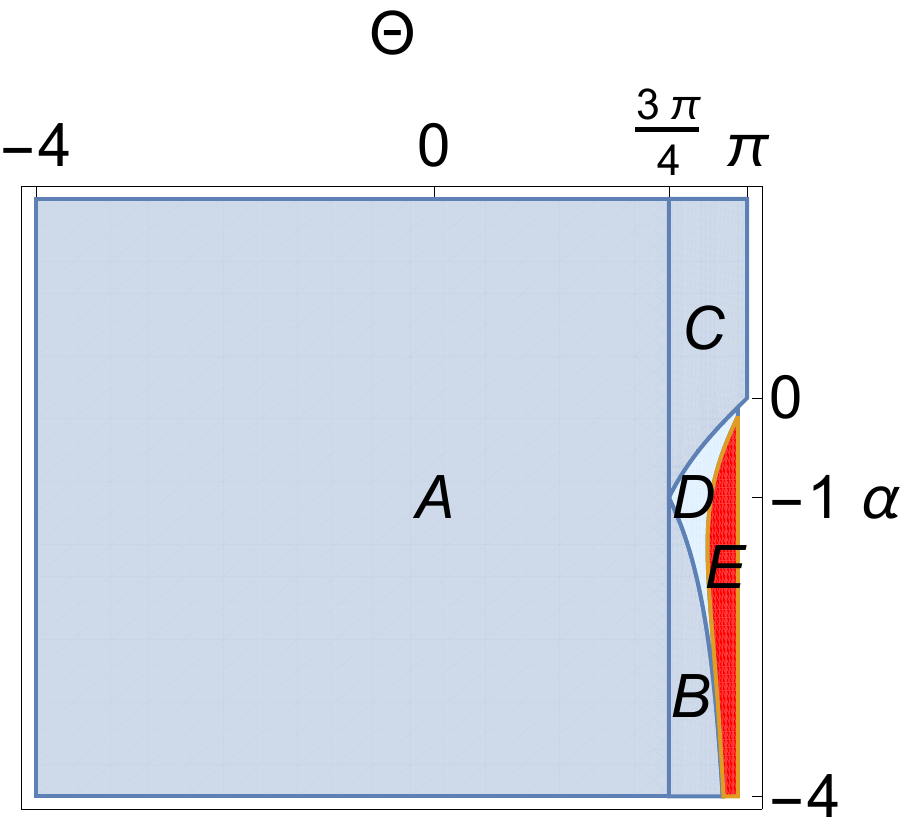}
\caption{
\label{fig:angular} 
Consider $K=-1,0,+1$. The second Robin eigenfunction of the geodesic disk $D_K(\Theta)$ with Robin parameter $\alpha$ has angular dependence in parameter regions $A,B,C,D$, by \autoref{basic}. The left side corresponds to hyperbolic disks ($K=-1$) and the right side to spherical caps ($K=+1$). For spherical caps with aperture $3\pi/4<\Theta<\pi$, region $B$ is where $\alpha \leq \cot \Theta$ and region $C$ is where $\alpha \geq \tan \Theta$. Angular dependence continues to hold in parameter region $D$ (see the comments on numerical work in \autoref{sec:construction}). In the exceptional region $E$, which corresponds to large spherical caps with suitably negative Robin parameter, the second Robin eigenfunction is not angular but instead (according to our numerical work) is radial.}
\end{center}
\end{figure}

We do not understand rigorously what happens in the spherical case when $3\pi/4<\Theta<\pi$ and $\alpha \in (\cot \Theta,\tan \Theta)$. According to numerical work, the second eigenfunction has angular dependence for a subset of that parameter region but the eigenfunction is instead radial for some parameter values, specifically when the cap fills almost the whole sphere ($\Theta$ close to $\pi$) and the Robin parameter lies in a certain negative range. \autoref{fig:angular} illustrates the proposition and our numerical findings. 

Regions $D$ and $E$ in the figure are found by numerically computing the lowest angular mode and lowest two radial modes on a spherical cap in order to determine where the second eigenfunction has angular dependence and where it is radial. The red region labeled ``$E$'' in the figure is the exceptional parameter set where the second eigenfunction is radial. It yields the red region in \autoref{fig:BS-FL}, after transforming the horizontal and vertical parameters to $t=4\pi \sin^2 \Theta/2$ and $\beta=(2\pi \sin \Theta) \alpha$. 
%

%See Bandle \cite[pp.\,122-124]{B80}, Ashbaugh and Benguria \cite[Section 6]{AB95}, or \autoref{basicproofs} below. 
A different proof for the Euclidean part (b) of the proposition was given by Freitas and Laugesen \cite{FL20}, using Bessel functions.  

\subsection*{Monotonicity in the radial direction for the first angular eigenfunction} 
\begin{figure}
\begin{center}
\includegraphics[scale=0.7]{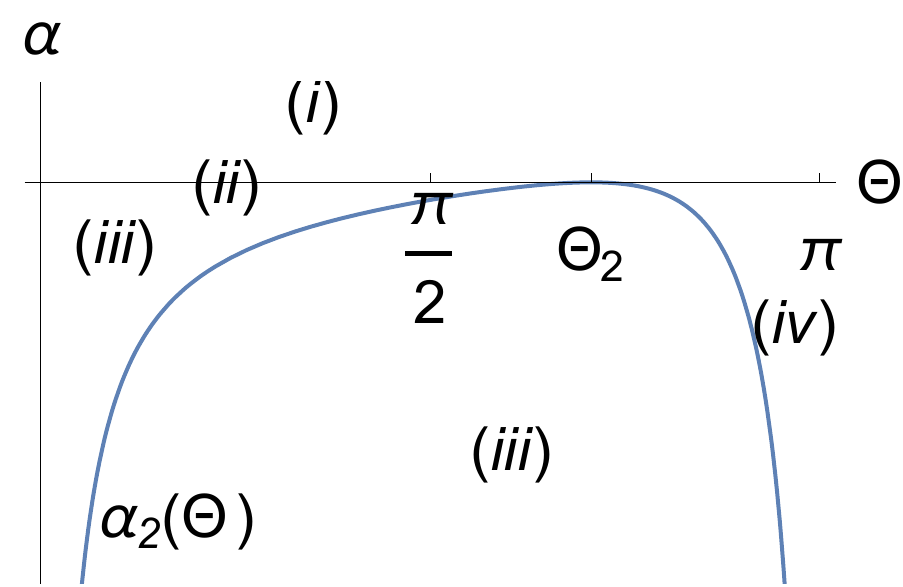}
\caption{
\label{fig:alpha2-plotlabeled} 
The regions for \autoref{basicsphere} parts (i)--(iv). 
}
\end{center}
\end{figure}

Next we aim for monotonicity properties of the radial part $g$ of the lowest eigenfunction having angular dependence. This eigenfunction has the form $u=g(\theta) \cos \phi$ or $g(\theta) \sin \phi$, since functions of the form $g(\theta) \cos m\phi, m \geq 2$, would generate larger eigenvalues; see Step 5 in the proof of \autoref{basic}. Importantly, the next results do not assume that this lowest angular eigenfunction gives the second eigenfunction.  
\begin{proposition}[Monotonicity of first angular eigenfunction: hyperbolic/Euclidean] \label{basichyperbolic} Suppose $K=-1$ or $K=0$, and let $\Theta>0$ and $\alpha \in \R$. 

If $g(\theta) \cos \phi$ is a first angular eigenfunction of $\Delta_K$ on a geodesic disk $D_K(\Theta)$ with Robin parameter $\alpha$, then one may take $g$ to be positive: $0=g(0)<g(\theta)$ whenever $\theta \in (0,\Theta]$. Furthermore: 

If $\alpha \leq 0$ then $g(\theta)$ is strictly increasing for $\theta \in (0,\Theta)$, with $g'>0$ there. 

If $\alpha>0$ then $g(\theta)$ first strictly increases and then strictly decreases: a maximum point $\theta_{max}$ exists such that $g'>0$ on $(0,\theta_{max})$ and $g'<0$ on $(\theta_{max},\Theta]$. 
\end{proposition}
The spherical case exhibits more complicated behavior, when the Robin parameter $\alpha$ is negative in part (iii) of the next proposition. Recall the aperture $\Theta_2 \simeq 0.70\pi$ and the function $\alpha_2(\Theta)$ that were defined before \autoref{th:bs}. Again we study the first angular mode, which is not necessarily the second eigenfunction. 
\begin{proposition}[Monotonicity of first angular eigenfunction: spherical] \label{basicsphere} Let $K=+1$ and take $\Theta \in (0,\pi)$ and $\alpha \in \R$. 

If $g(\theta) \cos \phi$ is a first angular eigenfunction of $\Delta_{+1}$ on a spherical cap of aperture $\Theta$ with Robin parameter $\alpha$, then one may take $g$ to be positive: $0=g(0)<g(\theta)$ whenever $\theta \in (0,\Theta]$. 

Furthermore, the behavior of $g$ depends on the sign of $\alpha$ as follows: 

\noindent (i) If $\alpha>0$ then $g(\theta)$ first strictly increases and then strictly decreases: a maximum point $\theta_{max}$ exists such that $g'>0$ on $(0,\theta_{max})$ and $g'<0$ on $(\theta_{max},\Theta]$. 

\noindent (ii) If $\alpha = 0$ and $0<\Theta \leq \Theta_2$ then $g(\theta)$ is strictly increasing, with $g'>0$ on $(0,\Theta)$. If $\alpha = 0$ and $\Theta_2<\Theta < \pi$ then $g(\theta)$ first strictly increases and then strictly decreases: a maximum point $\theta_{max}$ exists such that $g'>0$ on $(0,\theta_{max})$ and $g'<0$ on $(\theta_{max},\Theta)$. 

\noindent (iii) If $\alpha_2(\Theta) < \alpha < 0$ and $0<\Theta < \Theta_2$, or if $\alpha \leq \alpha_2(\Theta)$ and $0<\Theta <\pi$, then $g(\theta)$ is strictly increasing for $\theta \in (0,\Theta]$ with $g'>0$ on that interval (except $g^\prime$ vanishes at $\theta=\Theta_2$ when $\alpha = \alpha_2(\Theta)$ and $\Theta_2 \leq \Theta <\pi$). 

\noindent (iv) If $\alpha_2(\Theta) < \alpha < 0$ and $\Theta_2<\Theta <\pi$, then $g(\theta)$ first strictly increases and then strictly decreases and then strictly increases again: a local maximum point $\theta_{max}$ and local minimum point $\theta_{min}$ exist such that $g'>0$ on $(0,\theta_{max})$ and $g'<0$ on $(\theta_{max},\theta_{min})$ and $g'>0$ on $(\theta_{min},\Theta]$. 
\end{proposition}
\autoref{fig:alpha2-plotlabeled} illustrates the regions in the proposition. 

\subsection*{Relevant literature on the form of the eigenfunction} The Neumann and Dirichlet cases of the preceding propositions are known in all three constant curvature situations, by work of Bandle \cite[pp.\,122-124]{B80}, Ashbaugh and Benguria \cite[Section 3]{AB95}, Ashbaugh and Benguria \cite[p.~562]{AB93}, \cite[Section 3]{AB01}, Benguria and Linde \cite[Section 3]{BL07}. See the summary by Langford and Laugesen \cite{LL22a}, who completed the Neumann case by handling spherical caps larger than a hemisphere. 

The Robin case in curvature zero (\autoref{basic} for disks in Euclidean space) was treated by Freitas and Laugesen \cite[Section 5]{FL20}, \cite[Section 5]{FL21}, using explicit formulas for Bessel functions. For geodesic disks in hyperbolic space with $\alpha \in [-\sigma_1(\Theta),0)$, see Li, Wang and Wu \cite[Propositions 3.1 and 3.2]{LWW22}; here $\sigma_1(\Theta)$ is the first positive Steklov eigenvalue.

For spherical caps, we know of no prior work identifying properties of the second Robin eigenfunction or of the first angular Robin eigenfunction. 

The proofs below avoid special functions and instead rely on qualitative properties determined by the eigenfunction equation. 

\subsection*{Proof of \autoref{basic}} 
The first Robin eigenfunction is positive and hence by separation of variables it must be radial. Suppose $f(\theta)$ is a radial eigenfunction on the geodesic disk $D_K(\Theta)$ that is not the first eigenfunction, so that $f$ satisfies
\[
- \Delta_K f = \rho f , \qquad f^\prime(\Theta)+\alpha f(\Theta)=0 ,
\]
for some eigenvalue $\rho$, and $f$ changes sign since it is $L^2$-orthogonal to the first eigenfunction. The first four steps of this proof will show $\rho>\lambda_2(\Theta,\alpha)$ under the hypotheses of the proposition, so that the second Robin eigenfunction is definitely not radial.

\smallskip
Step 1. Observe $f^\prime(0)=0$ since the radial eigenfunction $f$ is smooth at the origin. Let
\[
v = f^\prime(\theta) \cos \phi 
\]
and notice $v \not \equiv 0$ since $f$ is nonconstant. This $v$ satisfies the eigenfunction equation $- \Delta_K v = \rho v$ with eigenvalue $\rho$, because 
\[
- \Delta_K v = \frac{d\ }{d\theta} (-\Delta_K f) \cos \phi = \rho f^\prime(\theta) \cos \phi = \rho v ,
\]
where the first equality relies on direct calculation and the Pythagorean identity $(\snK^{\! \prime})^2-\snK \snK^{\! \prime \prime}=1$.

\smallskip
Step 2. Suppose $f^\prime(\widetilde{\theta})=0$ for some $\widetilde{\theta} \in (0,\Theta]$, so that $v$ satisfies a Dirichlet condition on the boundary of the disk $D_K(\widetilde{\theta})$. Because $v$ changes sign (due to the factor $\cos \phi$), it cannot be the first Dirichlet eigenfunction of $\Delta_K$ on that disk and so $\rho$ must be a second or higher Dirichlet eigenvalue there. Hence by domain monotonicity for Dirichlet eigenvalues,  
\begin{align*}
\rho \geq \lambda_2(\widetilde{\theta},\infty) 
& \geq \lambda_2(\Theta,\infty) \\
& > \lambda_2(\Theta,\alpha) ,
\end{align*}
where the final inequality relies on strict monotonicity of the spectrum with respect to the Robin parameter. Thus $\rho > \lambda_2(\Theta,\alpha)$, as desired.

\smallskip
Step 3. Suppose next that $f^\prime \neq 0$ on $(0,\Theta]$, which means $f^\prime$ does not change sign. We may take $f^\prime>0$, so that the sign-changing property of $f$ implies $f(0)<0<f(\Theta)$. The Robin condition therefore implies
\[
\alpha = - \frac{f^\prime(\Theta)}{f(\Theta)} < 0 .
\]
Further, since $f(\theta_0)=0$ for some $\theta_0 \in (0,\Theta)$ we know $f$ is a Dirichlet eigenfunction on the disk $D_K(\theta_0)$ and so its eigenvalue must be positive: 
\[
\rho>0 .
\]

Let us determine the Robin condition satisfied by $v$. The eigenfunction equation $\Delta_K f + \rho f=0$ gives that 
\[
f^{\prime \prime} + \left( \frac{\snK^{\! \prime}}{\snK} + \rho \frac{f}{f^\prime} \right) f^\prime = 0.
\]
Evaluating at the boundary and using the Robin condition for $f$ shows that 
\[
f^{\prime\prime}(\Theta) + \gamma f^\prime(\Theta) = 0 
\]
where the constant is 
\[
\gamma = \ctK \Theta - \frac{\rho}{\alpha} 
\]
and we defined 
\[
\ctK \theta  = \frac{\snK^{\! \prime} \theta}{\snK \theta} = 
\begin{cases}
\cot \theta & \text{if $K=+1$,} \\
\theta^{-1} & \text{if $K=0$,} \\
\coth \theta & \text{if $K=-1$.}
\end{cases}
\]
Thus $v$ is a sign-changing Robin eigenfunction on $D_K(\Theta)$ with parameter $\gamma$ and eigenvalue $\rho$. It follows that 
\[
\rho \geq \lambda_2(\Theta,\gamma) .
\]
We want to show $\gamma>\alpha$, because then $\rho > \lambda_2(\Theta,\alpha)$. 

In the hyperbolic and Euclidean cases we have $\ctK \Theta \geq 0$. The same holds in the spherical case when $0<\Theta \leq \pi/2$. Thus in these cases, the proof that the second Robin eigenfunction is nonradial is complete, because $\gamma \geq -\rho/\alpha >0>\alpha$. 

\smallskip
Step 4. Consider now the spherical case ($K=+1$) with $\pi/2<\Theta<\pi$. If $\alpha \leq \cot \Theta$ then since $\gamma > \cot \Theta$ we have $\gamma > \alpha$, as needed. 

If $\tan \Theta \leq \alpha < 0$ then $\rho \geq 2$, as follows. The radial eigenfunction $f$ satisfies $-\Delta_{+1} f = \rho f$ with eigenvalue $\rho$ and Robin parameter $\alpha$ at $\theta=\Theta$, and the radial function $c(\theta) = - \cos \theta$ satisfies the eigenfunction equation $-\Delta_{+1} c = 2 c$ with eigenvalue $2$ and its Robin parameter at $\theta=\Theta$ is 
\[
- \frac{c^\prime(\Theta)}{c(\Theta)}  = \tan \Theta .
\] 
This Robin parameter is less than or equal to $\alpha$ by assumption and so the eigenvalue $2$ of $c$ is less than or equal to the eigenvalue $\rho$ of $f$, as we now justify. 

Suppose first that $f$ has its zero at some radius $\theta \in [\pi/2,\Theta)$. Then on the annulus between $\theta$ and $\Theta$ the function $f$ is a positive eigenfunction satisfying a Dirichlet condition (Robin parameter $+\infty$) at the inner boundary and a Robin condition at radius $\Theta$ with parameter $\alpha$, while on the same annulus, $c$ is a positive eigenfunction whose Robin parameter at $\Theta$ is less than or equal to $\alpha$; since positive eigenfunctions are automatically ground states, monotonicity of the spectrum with respect to the Robin parameter on each boundary portion implies that the eigenvalue of $c$ is less than or equal to that of $f$, meaning $2 \leq \rho$. 

Suppose next that $f$ has its zero at some radius $\theta \in (0,\pi/2)$. Then on the disk of radius $\theta$, the function $f$ is a negative eigenfunction satisfying a Dirichlet condition while $c$ is a negative radial eigenfunction satisfying some Robin condition at the boundary; hence again the eigenvalue of $c$ is less than or equal to that of $f$, giving $2 \leq \rho$ in this case too. 

By our assumption that $\alpha \geq \tan \Theta$ and the fact that $\rho \geq 2$, we obtain that  
\[
\gamma = \cot \Theta - \frac{\rho}{\alpha} \geq \frac{1}{\alpha} - \frac{2}{\alpha} > 0 > \alpha .
\]
Thus again $\gamma>\alpha$, as we wanted. 

Lastly, the assumption in \autoref{basic}(c) that $\alpha \notin (\cot \Theta,\tan \Theta)$ is needed only when $3\pi/4 < \Theta < \pi$, because if $\pi/2<\Theta \leq 3\pi/4$ then $\tan \Theta \leq \cot \Theta$ and the excluded interval is empty. 

\smallskip
Step 5. By a standard argument with separation of variables in the Rayleigh quotient, one finds that the second eigenfunction is a linear combination of some functions $g(\theta)\cos \phi$ and $g(\theta)\sin \phi$. (Angular factors $\cos m\phi$ and $\sin m\phi$ with $m \geq 2$ would give larger eigenvalues.) The Robin condition then says $g^\prime(\Theta)+\alpha g(\Theta)=0$.

\subsection*{Proof of \autoref{basichyperbolic} (Hyperbolic/Euclidean)}
Let $K=-1$ or $0$. The proposition concerns the first angular eigenfunction, which has the form of a radial function $g(\theta)$ multiplied by the angular part $\cos \phi$ or $\sin \phi$. Continuity of the eigenfunction at the origin demands that $g(0)=0$. Write $\lambda_{ang}(\Theta,\alpha)$ for this first angular eigenvalue. 

We begin by showing that after multiplying by $-1$ if necessary, one must have $g(\theta)>0$ when $0 < \theta \leq \Theta$. For suppose $g(\theta_0)=0$ for some $\theta_0 \in (0,\Theta]$. Then $v=g(\theta)\cos\phi$ is a Dirichlet eigenfunction with angular dependence on the disk $D_K(\theta_0)$, having eigenvalue $\lambda_{ang}(\Theta,\alpha)$. Hence 
\[
\lambda_{ang}(\Theta,\alpha) \geq \lambda_{ang}(\Theta,\infty) > \lambda_{ang}(\Theta,\alpha) ,
\]
where the first inequality holds by domain monotonicity of Dirichlet eigenvalues as we enlarge the disk $D_K(\theta_0)$ to $D_K(\Theta)$, and the second inequality holds by monotonicity of the eigenvalue with respect to the Robin parameter. This contradiction shows that $g(\theta_0) \neq 0$. Since $\theta_0$ was arbitrary, we see $g$ vanishes only at the origin, and so after replacing $g$ with $-g$ if necessary, we obtain that $g(\theta)>0$ whenever $0 < \theta \leq \Theta$. 

Next, applying the eigenfunction equation $-\Delta_K u = \lambda_{ang}(\Theta,\alpha) u$ to the eigenfunction $u=g(\theta)\cos \phi$ gives 
\begin{equation*} \label{eq:geqhyperbolic}
-\frac{1}{\snK \theta} \big( (\snK \theta) \, g^\prime \big)^{\! \prime} + \frac{1}{(\snK \theta)^2} g = \lambda_{ang}(\Theta,\alpha) g.
\end{equation*}
This equation holds for all $\theta>0$, since the ordinary differential equation is linear and so its solution extends to all positive $\theta$. 

Changing variable with
\[
s = 
\begin{cases}
\log (\tanh \theta/2) & \text{if $K=-1$} , \\
\log \theta & \text{if $K=0$} , 
\end{cases}
\]
we find that 
\begin{equation*} \label{eq:glogr}
\frac{d^2 g}{ds^2} = q g 
\end{equation*}
for $-\infty<s<\infty$ where  
\[
q(\theta) = 1-\lambda_{ang}(\Theta,\alpha) (\snK \theta)^2 .
\]
Notice $q$ is positive for small $\theta$ and so $g$ is a strictly convex function of $s$ near $-\infty$. Further, $g \to g(0)=0$ as $s \to -\infty$ and so $g$ must be increasing when $s$ is near $-\infty$. 

If $\lambda_{ang}(\Theta,\alpha) \leq 0$, then $q(\theta) \geq 1$ for all $\theta$. In particular, $g$ remains a strictly convex, strictly increasing function of $s$ all the way to the boundary, so that $g^\prime(\theta)>0$ on $(0,\Theta]$. Note that $g^\prime(\Theta)>0$ implies $\alpha<0$.

If $\lambda_{ang}(\Theta,\alpha) > 0$ then $q(\theta)$ is positive until $\theta$ becomes large enough that $q(\theta)$ changes sign and is thereafter negative. Thus $g$ is a strictly convex function of $s$ until it changes to become strictly concave, after which $g$ continues to be strictly concave for as long as it is positive. Thus either $dg/ds$ stays positive for the whole interval $\theta \in (0,\Theta)$ or else $dg/ds$ is first positive and then changes sign to remain negative through to the endpoint $\theta=\Theta$. That is, either $g'(\theta)>0$ for $\theta \in (0,\Theta)$, or else $g'>0$ on $(0,\theta_{max})$ and $g'<0$ on $(\theta_{max},\Theta]$. The first case has $g^\prime(\Theta) \geq 0$ and so $\alpha \leq 0$, while the second case has $g^\prime(\Theta)<0$ and hence $\alpha>0$. 

\subsection*{Proof of \autoref{basicsphere} (Spherical)} The proof that $g(0)=0$ and $g$ is positive on $(0,\Theta]$ proceeds exactly as for the hyperbolic/Euclidean case in the proof of \autoref{basichyperbolic}, and adapting that proof shows that $d^2 g/ds^2=qg$ where now 
\[
s = \log (\tan \theta/2) \in (-\infty,\infty) 
\]
and
\[
q(\theta) = 1-\lambda_{ang}(\Theta,\alpha) \sin^2 \theta .
\]
Again $q$ is positive for small $\theta$ and so $g$ is a strictly convex function of $s$ near $-\infty$, with $g \to g(0)=0$ as $s \to -\infty$, and so $g$ must be increasing when $s$ is near $-\infty$. 

\smallskip
\autoref{basicsphere} parts (i) and (ii). If $\alpha \geq 0$ then the Robin boundary condition forces $dg/ds$ to be nonpositive at the right endpoint $s(\Theta)$, and so at some point $g$ must switch from convex to concave. That is, $q(\theta)$ must change sign at least once on the interval $(0,\Theta)$. Noting that $\sin^2 \theta$ increases from zero before decreasing again to zero on $(0,\pi)$, we deduce $q$ has two roots $\theta_\pm$ satisfying $0<\theta_- < \pi/2 < \theta_+<\pi$ and that the smaller root $\theta_-$ must lie in $(0,\Theta)$. Write $s_\pm = s(\theta_\pm)$ for the $s$-values corresponding to the roots, so that 
\[
-\infty < s_- < 0 < s_+ < \infty
\]
and $s_- < s(\Theta) = \log \tan (\Theta/2)$. 

Suppose first that $\lambda_{ang}(\Theta,\alpha) \sin^2 \Theta \geq 1$, so that $q(\Theta) \leq 0$. Hence $q$ is positive on $(0,\theta_-)$ and negative on $(\theta_-,\Theta)$. The preceding paragraph shows that $g$ is strictly convex as a function of $s \in (-\infty,s_-)$ and strictly concave for $s \in (s_-,s(\Theta))$. If $\alpha>0$ then $dg/ds < 0$ at the endpoint $s(\Theta)$ and so we deduce that $g$ reaches a maximum at some point $\theta_{max}$ such that $g'>0$ on $(0,\theta_{max})$ and $g'<0$ on $(\theta_{max},\Theta]$. If $\alpha=0$ then $dg/ds = 0$ at the endpoint $s(\Theta)$ and so we deduce that $g$ is strictly increasing, with $g'>0$ on $(0,\Theta)$. 

Suppose next that $\lambda_{ang}(\Theta,\alpha) \sin^2 \Theta < 1$, so that $q(\Theta)>0$ and so $\theta_+ < \Theta$, with $q$ being positive on $(0,\theta_-)$, negative on $(\theta_-,\theta_+)$ and positive on $(\theta_+,\Theta)$. Our work above implies that $g$ is strictly convex as a function of $s$ on $(-\infty,s_-)$, strictly concave on $(s_-,s_+)$ and strictly convex on $(s_+,s(\Theta))$. Recalling that the slope $dg/ds$ is positive when $s$ is near $-\infty$ and is nonpositive at $s=s(\Theta)$, we deduce that for some number $s_{max} \in (-\infty,s(\Theta))$ one has $dg/ds>0$ on $(-\infty,s_{max})$ and $dg/ds<0$ on $(s_{max},s(\Theta))$. Determining $\theta_{max}$ from the relation $s_{max}=s(\theta_{max})$, we see $g^\prime(\theta)$ is positive on $(0,\theta_{max})$ and negative on $(\theta_{max},\Theta)$. 

Parts (i) and (ii) are now proved, noting for part (ii) in the proof above that the critical aperture $\Theta_2$ is defined so that $\lambda_{ang}(\Theta_2,0) \sin^2 \Theta_2 = 1$, with $\lambda_{ang}(\Theta,0) \sin^2 \Theta > 1$ when $\Theta \in (0,\Theta_2)$ and $\lambda_{ang}(\Theta,0) \sin^2 \Theta < 1$ when $\Theta \in (\Theta_2,\pi)$, using here \cite[Proposition 3.1]{LL22a} and the fact that $\lambda_{ang}=\lambda_2$ in the Neumann case $\alpha=0$ (by \autoref{basic}). 

\smallskip
\autoref{basicsphere} parts (iii) and (iv). Suppose $\alpha<0$ and $0<\Theta<\pi$. Note $g^\prime(\Theta)>0$ since the Robin parameter $\alpha$ is negative. 

Recall from \autoref{sec:regions} that the second Neumann eigenfunction $g_2(\theta) \cos \phi$ on the cap of aperture $\Theta_2$ is also a Robin eigenfunction on the cap of aperture $\Theta \in (0,\pi)$, with Robin parameter $\alpha_2(\Theta)$, and that $g_2^\prime>0$ except at $\Theta_2$ where $g_2^\prime$ vanishes. 

First suppose $\alpha < \alpha_2(\Theta)$, so that $g_2^\prime/g_2 < g^\prime/g$ at $\theta=\Theta$. \autoref{le:eigenineq} says that the same inequality must hold for all $\theta \in (0,\Theta]$, and so in particular $0 < g^\prime(\theta)$. Next suppose $\alpha = \alpha_2(\Theta)<0$. Then $g_2^\prime/g_2=g^\prime/g$ at $\theta=\Theta$ and so the equiality statement in the lemma implies that $g$ must be a positive multiple of $g_2$ and so $g^\prime(\theta)>0$ for all $\theta$ except $\theta=\Theta_2$. 

Now suppose $\alpha_2(\Theta) < \alpha < 0$ and $0<\Theta<\Theta_2$. Further suppose $g^\prime$ vanishes at some $\theta<\Theta$. It follows that $g$ must change from convex (as a function of $s$) to concave in order for $g^\prime$ to vanish, and then must change again to convex in order for $g^\prime$ to become positive at the endpoint $\Theta$. Since $q$ can change sign at most twice, we deduce $g$ must remain convex on $[\Theta,\pi)$ and hence also positive and strictly increasing there. In particular, $g^\prime(\Theta_2)>0$. The assumption $\alpha_2(\Theta) < \alpha$ means that $g_2^\prime/g_2 > g^\prime/g$ at $\Theta$, and so \autoref{le:eigenineq} implies that the same inequality must hold at $\Theta_2$, giving the contradiction  $g_2^\prime(\Theta_2)>0$. Therefore $g^\prime$ cannot vanish as we supposed, and hence $g^\prime>0$ on the whole interval $(0,\Theta]$. 

Finally, suppose $\alpha_2(\Theta) < \alpha < 0$ and $\Theta_2<\Theta<\pi$. The inequality $\alpha_2(\Theta) < \alpha$ means that $g_2^\prime/g_2 > g^\prime/g$ at $\Theta$, which by \autoref{le:eigenineq} implies the same inequality at $\Theta_2$, giving $g^\prime(\Theta_2) < 0$. Hence $g$ must change from strictly convex (as a function of $s$) to strictly concave in order for $g^\prime$ to be negative at $\Theta_2$, after which $g$ must change back to strictly convex in order to ensure $g^\prime$ is positive at the endpoint $\Theta$. It follows easily now that $g(\theta)$ first strictly increases and then strictly decreases and then strictly increases again, as claimed in part (iv) of the proposition. 

\subsection*{Properties of the radial part for aperture $\Theta_2$}
As above, $g_2$ is the radial part of the second Neumann eigenfunction for the spherical cap of aperture $\Theta_2$, and $g_2$ extends to a positive, increasing function on the whole interval $(0,\pi)$, as graphed in \autoref{fig:g-plot}. The next section needs the following facts about $\alpha_2(\Theta)=-g_2^\prime(\Theta)/g_2(\Theta)$.
\begin{lemma} \label{le:alpha2}
For $\Theta \in [\Theta_2,\pi)$ one has $\alpha_2(\Theta) > \cot \Theta$. Hence when $t \in [t_2,4\pi)$ one has $\beta_2(t) > 2\pi-t$, and also $\beta_2(t) \to -2\pi$ as $t \to 4\pi$. 
\end{lemma}
The lemma helps explain \autoref{fig:BS-FL}, where the graph of $\beta_2(t)$ lies above the upper boundary of region $\III$ and approaches height $-2\pi$ in the bottom right corner.
\begin{proof}
By applying the eigenfunction equation $\Delta_{+1}(g_2(\theta) \cos \phi) = -\lambda_2(\Theta_2,0) g_2(\theta) \cos \phi$ with Robin parameter $\alpha=0$ (for the Neumann boundary condition), one finds after some reorganization that 
\[
\frac{1}{\sin \theta} \Big( \frac{1}{\sin \theta} \, \big( g_2(\theta) \sin \theta \big)^\prime \Big)^{\! \prime} = -\lambda_2(\Theta_2,0) \frac{g_2(\theta)}{\sin \theta} < 0
\]
when $0<\theta<\pi$, and so 
\[
\frac{d^2\ }{dt^2} \big( g_2(\theta) \sin \theta \big) < 0
\]
where the new variable $t=t(\theta)=4\pi \sin^2 \theta/2$ satisfies $dt/d\theta=2\pi \sin \theta$. Thus $g_2(\theta) \sin \theta$ is a strictly concave function of $t$. Its first $\theta$-derivative at $\theta=\Theta_2$ equals 
\[
g_2^\prime(\Theta_2) \sin \Theta_2 + g_2(\Theta_2) \cos \Theta_2<0
\]
since $g_2^\prime(\Theta_2)=0$ by the Neumann boundary condition and $\Theta_2$ is larger than $\pi/2$ by definition in \autoref{sec:regions}. Hence the first $t$-derivative is negative at $t_2=t(\Theta_2)$, and so by concavity the $t$-derivative remains negative for all $t \in [t_2,4\pi)$, which means the first $\theta$-derivative is negative for all $\theta \in [\Theta_2,\pi)$. That is, 
\[
g_2^\prime(\theta) \sin \theta + g_2(\theta) \cos \theta < 0 ,
\]
which is equivalent to $\cot \theta<\alpha_2(\theta)$, as claimed in the lemma. 

Next, if $t \in [t_2,4\pi)$ then $t=4\pi \sin^2 \theta/2$ for some $\theta \in [\Theta_2,\pi)$. Multiplying the preceding inequality by $2\pi \sin \theta$ gives $(2\pi \sin \theta) \cot \theta < (2\pi \sin \theta) \alpha_2(\theta)$, which simplifies to $2\pi - t = 2\pi \cos \theta < \beta_2(t)$, giving the desired lower bound on $\beta_2$. Letting $t \to 4\pi$ implies $-2\pi \leq \liminf_{t \to 4\pi} \beta_2(t)$. 

To get an upper bound on the $\limsup$, let $0<\delta<\e<1$. By the proof of \autoref{basicsphere} above (choosing there the aperture $\Theta=\Theta_2$ and Robin parameter $\alpha=0$), we have $q(\theta) \to 1 > \e^2$ as $\theta \to \pi$ and hence $d^2 g_2/ds^2 > \e^2 g_2$ for all $\theta$ near $\pi$, where we recall that  $s = \log (\tan \theta/2)$ and $d/ds=(\sin \theta) d/d\theta$. Writing $y(s)=(1/g_2)dg_2/ds>0$, we deduce that 
\[
y^\prime(s) + y(s)^2 > \e^2 
\]
for all large $s$. Fix $s_0$ to be such a large number. If $0 < y(s_0) < \delta$ then the differential inequality for $y(s)$ forces its derivative to exceed a positive constant: $y^\prime(s) > \e^2 - \delta^2>0$, with this inequality holding not only at $s=s_0$ but on the whole open interval to the right of $s_0$ on which the value of $y(s)$ remains below $\delta$. Thus $y(s)$ must eventually exceed $\delta$ in value, after which its value remains above $\delta$, by invoking the differential inequality once more. Hence $y(s)>\delta$ for all large $s$, which means $(\sin \theta) (1/g_2) dg_2/d\theta > \delta$ for all $\theta$ near $\pi$ and hence $\beta_2(t) < -2\pi \delta$ for all $t$ near $4\pi$. Hence $\limsup_{t \to 4\pi} \beta_2(t) \leq -2\pi \delta$. Since $\delta<1$ was arbitrary, we have $\limsup_{t \to 4\pi} \beta_2(t) \leq -2\pi$ as desired, so that $\lim_{t \to 4\pi} \beta_2(t) = -2\pi$. 
\end{proof}

\section{\bf Proof of \autoref{th:bs} --- shape of the $\BS$ region}
\label{sec:bsproof}

Set $I$. The portion of set $I$ lying in the left halfplane is the semi-infinite strip 
\[
(-\infty,0) \times [-2\pi,0] = \{ (-4\pi \sinh^2 \Theta/2, \beta) : \Theta > 0, -2\pi \leq \beta \leq 0 \} ,
\]
where we have expressed $t$ as $4\pi \sinh^2 \Theta/2$. To satisfy the definition of the $\BS$ region, we show that \eqref{cond1} and \eqref{cond2a} hold in the hyperbolic case ($K=-1$) on a geodesic disk of radius $\Theta>0$ with Robin parameter $\alpha=\beta/(2\pi \sinh \Theta) \leq 0$. Indeed, for \eqref{cond1} the angular form $g(\theta) \cos \phi$ for the second eigenfunction holds by \autoref{basic}(a) while the positivity of $g$ and $g^\prime$ on $(0,\Theta)$ was shown in \autoref{basichyperbolic}.

The part of $I$ on the vertical axis is the interval $\{ (0,\beta) : -2\pi \leq \beta \leq 0 \}$. Let $\Theta > 0$ and note that \eqref{cond1} and \eqref{cond2a} hold in the Euclidean case ($K=0$) on a  disk of radius $\Theta$ with Robin parameter $\alpha=\beta/2\pi\Theta \leq 0$ by \autoref{basic}(b) and \autoref{basichyperbolic}. 

For the remainder of the proof we deal with sets in the open right halfplane. Take $K=+1$ from now on. The portion of set $I$ lying in the right halfplane can be expressed as the strip 
\[
\{ (4\pi \sin^2 \Theta/2, \beta) : 0<\Theta\leq\Theta_2, -2\pi \leq \beta \leq 0 \} 
\]
after writing $t=4\pi \sin^2 \Theta/2$. For such $\Theta$ and $\beta$ values, we see condition \eqref{cond1} holds by \autoref{basic}(c) since $\Theta_2<0.71\pi < 3\pi/4$ by \cite[Proposition 4.2]{LL22a}, while condition \eqref{cond2a} holds by \autoref{basicsphere}(ii)(iii).

\medskip
Set $\II$. This set lies in the right halfplane. Converting $t$ to $\Theta$, the set can be written as 
\[
\II = \{ (4\pi \sin^2 \Theta/2,\beta) : \Theta_2 < \Theta \leq 3\pi/4 , -2\pi \leq \beta \leq (2\pi \sin \Theta) \alpha_2(\Theta) \} .
\]
The angular form of the second eigenfunction for \eqref{cond1}, on a spherical cap of aperture $\Theta$ with Robin parameter $\alpha=\beta/(2\pi \sin \Theta) \leq \alpha_2(\Theta)$, holds by \autoref{basic}(c) since $\Theta \leq 3\pi/4$. Positivity of $g$ and of $g^\prime$ on $(0,\Theta)$ (except perhaps at one point, for $g^\prime$) follows from \autoref{basicsphere}(iii). 

\medskip
Set $\III$. This set may be written as 
\[
\III = \{ (4\pi \sin^2 \Theta/2,\beta) : 3\pi/4 < \Theta < \pi , -2\pi \leq \beta \leq 2\pi \cos \Theta \} 
\]
since $2\pi \cos \Theta = 2 \pi (1 - 2 \sin^2 \Theta/2) = 2 \pi - t$. The angular form of the second eigenfunction on a spherical cap of aperture $\Theta$ with Robin parameter $\alpha=\beta/(2\pi \sin \Theta) \leq \cot \Theta$ holds by \autoref{basic}(c). Positivity of $g$ and $g^\prime$ on $(0,\Theta)$ follows from \autoref{basicsphere}(iii) since  $\alpha \leq \cot \Theta < \alpha_2(\Theta)$ by \autoref{le:alpha2}. 

\medskip 
Set $\IV$. With $t=4\pi \sin^2 \Theta/2$, this final region in the right halfplane becomes 
\[
\IV = \{ (4\pi \sin^2 \Theta/2,\beta) : 3\pi/4 < \Theta < \pi , \tan \Theta \leq \beta/(2\pi \sin \Theta) \leq \alpha_2(\Theta) \} 
\]
where we used the definition of $\beta_2(t)$ in terms of $\alpha_2(\Theta)$ and used also that 
\begin{align*}
t \, \frac{4\pi-t}{2\pi-t} 
& = 4\pi \sin^2 \Theta/2 \frac{4\pi(1-\sin^2 \Theta/2)}{2\pi \cos \Theta} 
= 2\pi \sin \Theta \tan \Theta .
\end{align*}
By \autoref{basic}(c), the second eigenfunction has angular form as in \eqref{cond1} on a spherical cap of aperture $\Theta$ with Robin parameter $\alpha=\beta/(2\pi \sin \Theta) \geq \tan \Theta$. Positivity of $g$ and $g^\prime$ on $(0,\Theta)$ follows from \autoref{basicsphere}(iii) since  $\alpha \leq \alpha_2(\Theta)$, except that $g^\prime$ might vanish at one point. 

\medskip
Set $V$. The angular condition \eqref{cond1} holds by definition of set $V$, while the monotonic condition \eqref{cond2a} holds by \autoref{basicsphere}(iii) since set $V$ lies below the curve $\beta_2(t)$. 

\section{\bf Proof of \autoref{pr:FLdown} --- dropping downward in the $\FL$ set}
\label{sec:FLdownproof}

According to the angular hypothesis \eqref{cond1}, we may write $g$ and $g_*$, respectively, for the radial parts of the second eigenfunctions on the disk of radius $\Theta$ corresponding to Robin parameters $\alpha=\beta/(2\pi \, \snK \Theta)$ and $\alpha_*=\beta_*/(2\pi \, \snK \Theta)$. Both $g$ and $g_*$ are  positive on $(0,\Theta]$, by \autoref{basichyperbolic} and \autoref{basicsphere}. 

If the monotonicity condition \eqref{cond2a} holds then $g^\prime(\Theta) \geq 0$ and so $\beta \leq 0$ by the Robin boundary condition; recall also that $\beta \geq -2\pi$ by hypothesis in this proposition. Hence by the definition of the $\BS$ set, $(0,\beta) \in \BS$ if $K=0$ or $(4\pi \sin^2 \Theta/2,\beta) \in \BS$ if $K=+1$, so that the proposition is proved. Thus from now on we may suppose \eqref{cond2a} does not hold, so that $g^\prime$ changes sign and hence by \autoref{basicsphere}, the Up-Down-(Up) condition \eqref{cond2b} holds. 

To finish the proof, we must verify the remaining criterion for belonging to the set $\FL$, which is condition \eqref{cond3}. Normalize $g$ by a multiplicative constant so that the value at its local minimum point $\theta_{min}$ equals the value of $g_*$ at that point, meaning $g(\theta_{min})=g_*(\theta_{min})$. At the right endpoint $\theta=\Theta$, the Robin boundary condition and the assumption $\beta<\beta_*$ yield that 
\[
\frac{g_*^\prime(\Theta)}{g_*(\Theta)} = - \alpha_* < - \alpha = \frac{g^\prime(\Theta)}{g(\Theta)} .
\] 
Hence \autoref{le:eigenineq} implies $g_*^\prime/g_* < g^\prime/g$ for all $\theta \in (0,\Theta]$, so that $g_*^\prime < g^\prime$ at every point where $g_*=g$. In particular $g_*^\prime < g^\prime$ at $\theta_{min}$, and so by a short argument one concludes that $g_*>g$ on the interval $(0,\theta_{min})$.

Integrating the left side of \eqref{cond3} by parts yields that
\begin{align*}
& \int_0^{\theta_{min}} g(\theta) g^\prime(\theta) \left( \snK \frac{\theta}{2} \right)^{\! \! 2} \, d\theta \\
& = \frac{1}{4} \int_0^{\theta_{min}} ( g(\theta_{min})^2 - g(\theta)^2) \, \snK \theta \, d\theta \\
& > \frac{1}{4} \int_0^{\theta_{min}} ( g_*(\theta_{min})^2 - g_*(\theta)^2) \, \snK \theta \, d\theta \quad \text{since $g_*>g>0$} \\
& = \int_0^{\theta_{min}} g_*(\theta) g_*^\prime(\theta) \left( \snK \frac{\theta}{2} \right)^{\! \! 2} d\theta 
\end{align*}
by parts again. By the Up-Down-(Up) hypothesis \eqref{cond2b} for $g_*$ we know $g_*^\prime$ is positive until the local maximum ${\theta_*}_{max}$ of $g_*$, then negative until the local minimum ${\theta_*}_{min}$, and then positive again. Hence the last displayed integral is greater than or equal to $0$ if $\theta_{min} \leq {\theta_*}_{max}$, and if $\theta_{min} > {\theta_*}_{max}$ then the integral is greater than or equal to the integral over $(0,{\theta_*}_{min})$, which is nonnegative by the $\FL$ hypothesis \eqref{cond3} for $g_*$. Thus in either case the last displayed integral is nonnegative. Hence, as we needed to prove, condition \eqref{cond3} holds for $g$. 

\section{\bf Construction of \autoref{fig:BS-FL} and \autoref{fig:contourplot}}
\label{sec:construction}

Readers are encouraged to download a Mathematica notebook \cite{LL23notebook} in order to follow along with the explanations below of how the figures were created. 

\subsection*{Construction of \autoref{fig:BS-FL}} \textbf{$\BS$ set.} The sets $I$--$\IV$ appearing in \autoref{fig:BS-FL} are specified in \autoref{th:bs} and can be plotted straightforwardly in Mathematica. The curves $\beta_2$ and $\beta_4$ and the set $V$ require some explanation. 

The curve $\beta_2(t) = -(2\pi \sin \Theta) \, g_2^\prime(\Theta)/g_2(\Theta)$ that forms the upper boundary for sets $\II, \IV, V$ is defined prior to the statement of \autoref{th:bs}, with $t$ and $\Theta$ related by $t=4 \pi \sin^2 \Theta/2$ and with $g_2(\theta)=P_{n_2}^{-1}(\cos \theta)$ defined in terms of an associated Legendre function. The parameter $n_2$ for the Legendre function is known theoretically to exist and be the unique value such that $g_2$ increases until some aperture $\Theta_2$ at which $g_2^\prime=0$, and $g_2$ continues to increase thereafter. Thus the curve $\beta_2(t)$ is negative until $t_2$, where it equals $0$, after which $\beta_2(t)$ becomes negative again. An exact formula for $n_2$ is not known, but by numerical experimentation one finds an accurate approximation to be $n_2 \simeq 0.851187$, which is the value used to plot $\beta_2(t)$ in the figure. 

The curve $\beta_4(t)$ is constructed analogously to $\beta_2(t)$, except with $g_4(\theta)=P_{n_4}^{-1}(\cos \theta)$ where $n_4 = 0.908729$. This parameter choice is an approximation to the largest $n$ for which the Front-Loaded condition \eqref{cond3} holds with $\beta=0$; for explanation, see our Neumann result \cite[proof of Theorem 1.1]{LL22a}. The curve $\beta_4(t)$ crosses the horizontal axis at $t_4 \simeq 11.828$. Note our earlier Neumann result used the slightly smaller number $(16/17) 4\pi$ in order to obtain a rigorous result.  

To plot the set $V$ in \autoref{fig:BS-FL}, we needed to identify the parameter values $(t,\beta)$ at which the angular condition \eqref{cond1} holds, that is, at which the first eigenvalue among angular modes is smaller than the second eigenvalue among radial modes. These eigenvalues were computed numerically in Mathematica with the NDEigenvalues command, after substituting $g(\theta)=\rho(\theta) h(\theta)$ for a suitably chosen weight $\rho$ in order to convert the Robin condition on $g$ into a Neumann condition on $h$, at $\theta=\Theta$. 
 
\textbf{$\FL$ set.} To determine the set $\FL$, one must verify the angular condition \eqref{cond1}, Up-Down(-Up) condition \eqref{cond2b}, and Front-Loaded condition \eqref{cond3}. In the first quadrant, where $t \geq 0$ and $\beta > 0$, the angular condition holds by \autoref{basic} and the Up-Down condition with $\theta_{min}=\Theta$ holds by \autoref{basichyperbolic} and \autoref{basicsphere}(i). Thus only the Front-Loaded condition need be checked. To avoid numerical differentiation, first integrate by parts in \eqref{cond3}. Then given a $t$ value, evaluate the corresponding $\Theta$ and apply a numerical bisection method to find the largest $n$ value for which \eqref{cond3} holds for $g(\theta)=P_n^{-1}(\cos \theta)$ with $\theta_{min}=\Theta$. This $n$ value and $\Theta$ determine $\beta$ from the Robin boundary condition $\beta/(2\pi \sin \Theta) =-g^\prime(\Theta)/g(\Theta)$. By \autoref{pr:FLdown}(ii), the $\FL$ set contains the segment dropping down from the point $(t,\beta)$ to the horizontal axis, and so $(t,\beta)$ lies on the upper boundary of the $\FL$ set. Performing this procedure for 60 reasonably-spaced $t$-values between $0$ and $t_4$ and then joining the resulting points yields an accurate representation of the $\FL$ set in the first quadrant in \autoref{fig:BS-FL}. 

Next we handle the part of the $\FL$ set in the fourth quadrant. The region above the curve $\beta_2(t)$ and to the right of $t_2$ satisfies the Up-Down-Up condition by \autoref{basicsphere}(ii)(iv). Points on the curve $\beta_4(t)$ satisfy the Front-Loaded condition by the choice of $n_4$, and satisfy the angular condition too provided we stay above the graph of $t(4\pi-t)/(2\pi-t)$. Thus the region in the fourth quadrant bounded by that graph, the horizontal axis and the graphs of $\beta_2$ and $\beta_4$ belongs to $\FL$, by dropping downward from the graph of $\beta_4$ via \autoref{pr:FLdown}(ii). Finally, the same reasoning also applies below the graph of $t(4\pi-t)/(2\pi-t)$ provided the angular condition can be verified numerically, thus obtaining the additional small piece of the $\FL$ set in \autoref{fig:BS-FL}. 

\subsubsection*{Which aspects above required numerical work?} In the $\FL$ region, numerics were needed primarily to verify the Front-Loaded condition, while in the $\BS$ region for set $V$, numerics were needed to verify that the second mode is angular. 

\subsubsection*{Remark} 
The Front-Loaded condition \eqref{cond3} holds for $g_4$ at $(t_4,0)$, as mentioned above for the Neumann condition. Hence the Front-Loaded condition holds at $(t,\beta_4(t))$ for all $t<t_4$ (smaller apertures) because the integral in \eqref{cond3} would include less negative contribution than when $t=t_4$, and it holds for all $t_4<t<\pi$ (larger apertures) because the integral would include more positive contributions. Similar reasoning can be applied on other curves constructed like $\beta_4$ but with different values of $n$. By combining this approach with \autoref{pr:FLdown}(ii), one may justify large parts of the $\FL$ set by evaluating the integral in \eqref{cond3} at just finitely many points $(t,\beta)$ in the first quadrant. In this sense, the derivation of the $\FL$ set shown in \autoref{fig:BS-FL} can be regarded as partly numerical and partly rigorous.  

\subsection*{Construction of \autoref{fig:contourplot}: contour plot of the lowest angular eigenvalue} 
The first angular mode on a spherical cap has the form $u=g(\theta) \cos \phi$ (see \autoref{eigenstuff}) where $g(\theta)=P_n^{-1}(\cos \theta)$ (see \autoref{sec:Legendre}) with an analogous formula in the hyperbolic case. Recalling that $t=4\pi \sin^2 \Theta/2$, we see that each point $(t,\beta(t))$ along the graph of $\beta(t) = -(2\pi \sin \Theta) \, g^\prime(\Theta)/g(\Theta)$ corresponds to a spherical cap with aperture $\Theta \in (0,\pi)$ and Robin parameter $\alpha(\Theta)=-g^\prime(\Theta)/g(\Theta)$. The Robin eigenvalue is the same on each of these caps, namely $\lambda=n(n+1)$ by \autoref{sec:Legendre}, since the underlying eigenfunction $u$ is the same for each cap. Thus the graph $(t,\beta(t))$ is a level curve or contour for the eigenvalue of the first angular mode.

\section*{Acknowledgments}
Richard Laugesen's research was supported by grants from the Simons Foundation (\#964018) and the National Science Foundation ({\#}2246537).

%\section*{Conflict of interest statement }
%
%On behalf of all authors, the corresponding author states that there is no conflict of interest.
%
%\section*{Data availability statement}
%
%Data sharing is not applicable to this article as no datasets were generated or analyzed during the current study.

\appendix

\section{\bf Legendre functions --- radial and angular modes} \label{sec:Legendre}

Separation of variables yields eigenfunctions as stated below on the $2$-dimensional sphere and hyperbolic space, in terms of the radial variable $\theta$ and angular variable $\phi$. The eigenfunction equation $-\Delta_K u = \lambda u$ (see \autoref{sec:coordinates}) can be verified straightforwardly for the functions below, using the associated Legendre ODE for $y=P_n^{-m}$ (see \cite[eq.\,14.2.2]{DLMF}):  
\[
(1-x^2) y^{\prime \prime}(x) - 2x y^\prime(x) + \left( n(n+1) - \frac{m^2}{1-x^2} \right) y(x) = 0  .
\]

\subsection*{Spherical eigenfunctions ($K=+1$)} 
\begin{align*}
u & = P_n^{-m}(\cos \theta) e^{i m\phi} , \qquad m=0,1,2,\ldots , \\
\lambda & = n(n+1) = 
\begin{cases}
-\frac{1}{4} - k^2 & \text{when $n = -\frac{1}{2} +ik$ with $k \geq 0$,} \\
-\frac{1}{4} + k^2 & \text{when $n = -\frac{1}{2} +k$ with $k \geq 0$.}
\end{cases}
\end{align*}
\noindent For radial modes one takes $m=0$, while $m=1$ yields the first angular mode. The eigenvalue $\lambda$ is determined (implicitly) when a boundary condition is imposed on $u$. 

\subsection*{Hyperbolic eigenfunctions ($K=-1$)} 
\begin{align*}
u & = i^{-m} P_n^{-m}(\cosh \theta) e^{i m\phi} , \qquad m=0,1,2,\ldots , \\
\lambda & = -n(n+1) = 
\begin{cases}
\frac{1}{4} + k^2 & \text{when $n = -\frac{1}{2} +ik$ with $k \geq 0$,} \\
\frac{1}{4} - k^2 & \text{when $n = -\frac{1}{2} +k$ with $k \geq 0$.}
\end{cases}
\end{align*}

\subsection*{Steklov case $\lambda=0$} The Steklov spectrum on a surface consists of the (negatives of the) Robin parameters whose corresponding eigenvalues equal $0$. 

In the spherical case, making the choice $n=0$ (so that $\lambda=0$) yields 
\[
u = P_0^{-m}(\cos \theta) e^{i m\phi}=\frac{1}{m!} \left( \tan \frac{\theta}{2} \right)^{\! \! m} e^{im\phi} ,
\]
which is analogous to the usual Steklov eigenfunction $r^m e^{im\phi}$ in the Euclidean case. In particular, when $m=1$ the eigenfunction $(\tan \theta/2) e^{i\phi}$ has Robin parameter $\alpha=-(\tan \theta/2)^\prime/(\tan \theta/2)=-1/\sin \Theta$ at aperture $\theta=\Theta$.

Similarly, the hyperbolic case yields a Steklov eigenfunction
\[
u = i^{-m} P_0^{-m}(\cosh \theta) e^{i m\phi}=\frac{1}{m!} \left( \tanh \frac{\theta}{2} \right)^{\! \! m} e^{im\phi} ,
\]
and when $m=1$ the Robin parameter at geodesic radius $\Theta$ is $\alpha=-1/\sinh \Theta$.

\section{\bf A lemma relating eigenvalues and endpoint values} \label{sec:usefullemma}

The next lemma relates the eigenvalues to the endpoint values of the eigenfunctions, for the lowest ``angular'' mode. Recall the function $\snK \theta$ defined earlier in \eqref{eq:snKdef}. 
\begin{lemma} \label{le:eigenineq}
Fix $K = -1,0,+1$ and $\lambda,\lambda_* \in \R$. Assume $\Theta>0$, and further suppose when $K=+1$ that $0<\Theta<\pi$. Suppose $g,g_* \in C^2[0,\Theta]$ satisfy
\begin{align*}
- \frac{1}{\snK \theta} \big( (\snK \theta) g^\prime(\theta) \big)^\prime + \frac{1}{(\snK \theta)^2} \, g(\theta) & =  \lambda g(\theta) , 
\\
- \frac{1}{\snK \theta} \big( (\snK \theta) g_*^\prime(\theta) \big)^\prime + \frac{1}{(\snK \theta)^2} \, g_*(\theta) & =  \lambda_* g_*(\theta) , 
\end{align*}
when $\theta \in (0,\Theta)$. If $g$ and $g_*$ are positive on $(0,\Theta]$ then 
\[
\sign \left( \frac{g_*^\prime}{g_*} - \frac{g^\prime}{g} \right) \!  = \sign(\lambda-\lambda_*) \qquad \text{on $(0,\Theta]$} .
\]
Thus if $g_*^\prime/g_* < g^\prime/g$ for some $\theta \in (0,\Theta]$ then that inequality holds for all $\theta$, and if $g_*^\prime/g_* = g^\prime/g$ for some $\theta \in (0,\Theta]$ then equality holds for all $\theta$ and hence $g_*=(\text{const.})g$. 

\end{lemma}
Results of this kind are well known. We include a short proof for the reader's convenience. 
\begin{proof}
Multiply the differential equation for $g$ by $(\snK \theta) g_*(\theta)$ and multiply the differential equation for $g_*$ by $(\snK \theta) g(\theta)$, and then subtract and integrate from $0$ to $\tau$, for an arbitrary $\tau \leq \Theta$. Hence 
\[
\int_0^\tau \left( g(\theta) \big( (\snK \theta) g_*^\prime(\theta) \big)^{\! \prime} - g_*(\theta) \big( (\snK \theta) g^\prime(\theta) \big)^{\! \prime} \right) \! d\theta = \left( \lambda - \lambda_* \right) \int_0^\tau g(\theta) g_*(\theta) \snK \theta \, d\theta .
\]
On the right side, the integral is positive since $g$ and $g_*$ are positive, recalling also when $K=+1$ that $\tau \leq \Theta<\pi$. Thus the sign of the right side equals $\sign(\lambda - \lambda_*)$. 

The fundamental theorem evaluates the left side to 
\begin{equation*} \label{eq:fundthspherical}
g(\theta) (\snK \theta) g_*^\prime(\theta) - g_*(\theta) (\snK \theta) g^\prime(\theta) \Big|_0^\tau = (\snK \tau) g(\tau) g_*(\tau) \left( \frac{g_*^\prime}{g_*} - \frac{g^\prime}{g} \right) \! (\tau) .
\end{equation*}
Since $\snK \tau , g(\tau)$ and $g_*(\tau)$ are positive, the sign of this side equals $\sign (g_*^\prime/g_* - g^\prime/g)(\tau)$. The remaining statements in the lemma follow easily. 
\end{proof}

\bibliographystyle{plain}

\end{document}